\documentclass{tac} 
\usepackage{latexsym}
\usepackage{amsfonts}
\usepackage{amssymb}
\usepackage{amsmath}
\usepackage{amscd}
\usepackage{eucal}
\usepackage[all]{xy}
\usepackage{pstricks}
\newcommand{\bs}{\boldsymbol}
\newcommand{\mb}{\mathbf}
\renewcommand{\dot}{\centerdot}

\begin{document}


     {}
     {}


\title{Biequivalences in tricategories}
\author{Nick Gurski}
\address{Department of Pure Mathematics, University of Sheffield, Sheffield, UK, S3 7RH}
\eaddress{nick.gurski@sheffield.ac.uk}

\maketitle

\begin{abstract}
We show that every internal biequivalence in a tricategory $T$ is
part of a biadjoint biequivalence.  We give two applications of this result, one for transporting monoidal structures and one for equipping a monoidal bicategory with invertible objects with a coherent choice of those inverses.
\end{abstract}


\section*{Introduction}

It is common in mathematics to regard two objects $X$ and $Y$ as being the same if there is an isomorphism between them, that is a pair of maps $f:X \rightarrow Y$ and $g:Y \rightarrow X$ such that $fg = 1_{Y}$ and $gf = 1_{X}$.  More importantly, any specific isomorphism $f:X \rightarrow Y$ gives an explicit means for transporting information about $X$ to information about $Y$.  The choice of ambient category plays a very important role in this process, especially when we examine two categories with the same objects but different morphisms.  For example, first take the category of CW-complexes and continuous maps.  In this category, an isomorphism between $X$ and $Y$ gives a formula for transporting cell structures from $X$ to $Y$.  If we now take the category of CW-complexes and homotopy classes of continuous maps, then isomorphisms are now homotopy equivalences.  An isomorphism then no longer gives a recipe for transporting cell structures, but instead only gives a recipe for transporting homotopical information like homotopy and homology groups.

An important property of an isomorphism $f:X \rightarrow Y$ is that it uniquely determines the inverse $g:Y \rightarrow X$ by the formulas $fg = 1_{Y}$ and $gf = 1_{X}$; the proof is exactly the same as that showing that any group element has a unique inverse.  Thus we have two concepts which are \textit{a priori} different:
\begin{itemize}
\item the property that $f:X \rightarrow Y$ is an isomorphism, and
\item the structure consisting of a pair $(f,g)$ of morphisms $f:X \rightarrow Y$, $g:Y \rightarrow X$ such that $fg = 1_{Y}$ and $gf = 1_{X}$.
\end{itemize}
We have then uncovered that a morphism $f$ has the property of being an isomorphism if and only if there is a pair $(f,g)$ with the isomorphism structure, and that moreover the pair $(f,g)$ is uniquely determined by $f$ alone.  This is an example of a structure (being part of an isomorphism pair) that is determined in a unique way by a property (the existence of an inverse).

Moving up to the case in which the objects of study are now the 0-cells of some 2-category (or more generally, some bicategory), there are now more possible notions of sameness.  While we can ask if two 0-cells of a 2-category are isomorphic, it is much more common to ask if they are equivalent.  The canonical example of a 2-category is that of categories, functors, and natural transformations, and in this 2-category we see from experience that equivalence is the natural notion of sameness.

A functor $F:X \rightarrow Y$ is often defined to be an equivalence if it is essentially surjective, full, and faithful.  It is then shown that a functor $F$ is an equivalence if and only if there exists a functor $G:Y \rightarrow X$ such that the composites $FG, GF$ are naturally isomorphic to the identity functors $1_{Y}, 1_{X}$, respectively.  This definition identifies a property of a functor $F$ that allows us to conclude that two categories $X$ and $Y$ are, in some sense, the same.

On the other hand, we can make the definition of the structure of adjoint equivalence $F \dashv_{eq} G$ which consists of a functor $F:X \rightarrow Y$, a functor $G:Y \rightarrow X$, and two natural isomorphisms $\eta: 1_{X} \Rightarrow GF, \varepsilon: FG \Rightarrow 1_{Y}$ such that the following two diagrams commute.
\[
\xy
{\ar^{F \eta} (0,0)*+{F}; (25,0)*+{FGF} };
{\ar^{\varepsilon F} (25,0)*+{FGF}; (25,-15)*+{F} };
{\ar_{1} (0,0)*+{F}; (25,-15)*+{F} };
{\ar^{\eta G} (50,0)*+{G}; (75,0)*+{GFG} };
{\ar^{G \varepsilon} (75,0)*+{GFG}; (75,-15)*+{G} };
{\ar_{1} (50,0)*+{G}; (75,-15)*+{G} }
\endxy
\]
It is clear that if $F \dashv_{eq} G$ is an adjoint equivalence, then $F$ is an equivalence.  It is also well-known \cite{mac} that every equivalence $F$ can be completed to an adjoint equivalence.  Now $F$ no longer determines $G, \eta, \varepsilon$ uniquely, but instead only determines them up to a unique isomorphism preserving the adjoint equivalence structure.  We refer the reader to the paper \cite{kl} for a general discussion of this phenomenon in the setting of algebras for a 2-monad.

The aim of this paper is to establish an analogous result in three-dimensional category theory.  We define two different notions of sameness internal to a tricategory, one a property of a 1-cell and the other a structure involving 1-, 2-, and 3-cells satisfying certain axioms.  These two notions are that of a 1-cell having the property ``is a biequivalence'' on the one hand, and the structure of a biadjoint biequivalence on the other hand.  Our main result is that, in any tricategory $T$, every 1-cell which is a biequivalence is part of a biadjoint biequivalence.

The proof of this result proceeds in three steps.  First, we show that it is true for the special case when $T=\mb{Bicat}$, the tricategory of bicategories, functors, pseudo-natural transformations, and modifications.  Second, we prove a result about transporting biadjoint biequivalences; more precisely, we show that if a functor $F:S \rightarrow T$ satisfies a kind of local embedding condition and $T$ has the property that every biequivalence is part of a biadjoint biequivalence, then $S$ has that property as well.  Finally, we prove that this property is inherited by functor tricategories from the target, so that if $T$ is a tricategory in which every biequivalence is part of a biadjoint biequivalence then the same holds for the functor tricategory $\mb{Tricat}(S,T)$ for any $S$.  The main result then follows from these theorems and coherence for tricategories by considering the Yoneda embedding.

We also give two applications of this result.  The first is a transport-of-structure result, showing how biequivalences $F:X \rightarrow Y$ between bicategories can be used to transport monoidal structures from $X$ to $Y$.  This relies on choosing a weak inverse $G:Y \rightarrow X$ for the definition of the tensor product on $Y$, and then requires the rest of the biadjoint biequivalence structure in order to define the higher cells that are part of the definition of a monoidal bicategory and to check that they satisfy the necessary axioms.  We also indicate how to prove similar results when the monoidal structure is replaced with a braided monoidal, sylleptic monoidal, or symmetric monoidal one.

The second application of the main result is to an elucidation of those monoidal bicategories in which every object is weakly invertible, called Picard 2-categories here.  We show that every Picard 2-category is monoidally biequivalent to one in which a coherent choice of inverses has been made.  In fact, the result is much stronger in that we show that the forgetful functor from coherent Picard 2-categories (those with a choice of inverses) to Picard 2-categories (those monoidal bicategories which merely have the property that every object is invertible) is a triequivalence.  We go on to improve this result by defining functors that preserve a given coherent structure up to equivalence, and show that every monoidal functor between Picard 2-categories can be given the structure of such.

The paper is organized as follows.  The first section is a warm-up in which we give a proof of the fact that every equivalence in a bicategory is part of an adjoint equivalence.  We do this to give the reader a taste of the strategy that will be used later in the tricategorical case so as to clearly indicate the crucial points.  This material is well-known although I am unaware of a reference that presents this result in full detail using the argument we give below.

The second section gives the definitions of biequivalence and biadjoint biequivalence that are at the heart of this paper.  Both of these we express in the completely general case, working in an arbitrary tricategory $T$.  The definition of biadjoint biequivalence has two forms, with and without the ``horizontal cusp'' axioms, and we discuss briefly why these two definitions are logically equivalent using the calculus of mates.

Section 3 gives a proof of our main result in the special case where $T= \mb{Bicat}$.  This proof is largely calculation, much as the proof that every equivalence in $\mb{Cat}$ is part of an adjoint equivalence is done by straightforward calculation in Section 1.  We use that a biequivalence in $\mb{Bicat}$ can be characterized in two different ways:  as a functor having a weak inverse or as a functor which is biessentially surjective and a local equivalence.  The proof largely consists of using this alternate characterization of biequivalences, results from section 1, and the biadjoint biequivalence axioms to construct the other cells of the biadjoint biequivalence.  This is done by constructing the components on objects first, and then building up the rest of the structure afterwords.  The reader will note that the technical difficulties lie in two places: in constructing these cells once the components on objects are given, and in checking that our constructions satisfy all of the biadjoint biequivalence axioms since we only need a subset of them in order to define all the needed cells.

The fourth section provides the proof of our main result, that every biequivalence in a tricategory is part of a biadjoint biequivalence.  The proofs in this section rely very heavily on coherence for tricategories in the form ``every diagram in a free tricategory commutes'' to simplify the pasting diagrams required.  We often leave the particulars of checking axioms to the reader as the diagrams are very large, but we state exactly which axioms are needed in each case.  The bulk of the technical work goes into showing that the components we construct satisfy the axioms required to be transformations, modifications, or perturbations, while the biadjoint biequivalence axioms are immediate.

Sections 5 and 6 give the two applications mentioned above.  Both of these results should be thought of in the form ``a certain 3-dimensional monad has property P'' in each case.  In the case of lifting monoidal structures, the monad would be the free monoidal bicategory monad on the tricategory $\mb{Bicat}$, and the property P would be a 3-dimensional version of flexibility \cite{bkp}.  In the case of Picard 2-categories, the monad would be the free Picard 2-category monad on the tricategory of monoidal bicategories, and the property P would be a kind of 3-dimensional idempotency \cite{kl}.  This is the proper way to view these results, although we do not pursue the details here because of the lack of groundwork on 3-dimensional monads on tricategories.

The author would like to thank Bruce Bartlett for his interest in these results, and John Baez and Eugenia Cheng for enlightening discussions.

\section{Equivalences and adjoint equivalences in bicategories}

We begin by reviewing the relevant results for bicategories that
we will later generalize to tricategories.  We assume that the
reader is familiar with bicategories and the coherence theorem for
bicategories (see \cite{mp} for coherence for bicategories, or \cite{js} for a discussion of coherence, including functors, for the case of monoidal categories instead of bicategories).  We begin
with some basic definitions.

\definition
Let $B$ be a bicategory, and let $f: x \rightarrow y$ and $g:y
\rightarrow x$ be 1-cells in $B$.  An \textit{adjunction} $f
\dashv g$ consists of a 2-cell $\varepsilon: fg \Rightarrow 1_{y}$
and a 2-cell $\eta: 1_{x} \Rightarrow gf$ such that the following
two diagrams (the triangle identities) commute.
\begin{displaymath}
\xymatrix{
g \ar[r]^{l_{g}^{-1}} \ar[rrrrd]_{1_{g}} & 1_{x}g \ar[r]^{\eta * 1_{g}} & (gf)g \ar[r]^{a} & g(fg) \ar[r]^{1_{g}*\varepsilon} & g1_{y} \ar[d]^{r} \\
&&&&g}
\end{displaymath}
\begin{displaymath}
\xymatrix{
f \ar[r]^{r_{f}^{-1}} \ar[rrrrd]_{1_{f}} & f 1_{x} \ar[r]^{1_{f}*\eta} & f(gf) \ar[r]^{a^{-1}} & (fg)f \ar[r]^{\varepsilon * 1_{f}} & 1_{y}f \ar[d]^{l} \\
&&&&f}
\end{displaymath}
We then say that $f$ is left adjoint to $g$, or that $g$ is right
adjoint to $f$.
\enddefinition

\remark  In the bicategory $\mathbf{Cat}$, the associativity and unit isomorphisms are all identities.  In that case, this definition reduces to the usual definition of an adjunction between functors.
\endremark

\definition
An adjunction $f \dashv g$ is an \textit{adjoint equivalence} if
$\varepsilon$ and $\eta$ are invertible.  In this case,
we write $f \dashv_{eq} g$.
\enddefinition

\theorem \label{cat} Let $F: X \rightarrow Y$ and $G:Y \rightarrow X$ be functors, and let $\alpha:FG \Rightarrow 1_{Y}$ and $\beta: 1_{X} \Rightarrow GF$ be natural isomorphisms.  Then there is a unique adjoint equivalence $(F, G, \varepsilon, \eta)$ in $\mathbf{Cat}$ such that $\varepsilon = \alpha$.
\endtheorem
\begin{proof}
Let $\varepsilon = \alpha$.  The second triangle identity states
that $\varepsilon F \circ F\eta = 1_{F}$.  By the invertibility of
$\varepsilon$, this equation is the same as $F \eta = (\varepsilon
F)^{-1}$.  The righthand side of this equation is well-defined and
$F$ is full and faithful since it is an equivalence of categories,
so we define $\eta_{x}: x \rightarrow GFx$ to be the unique arrow
such that $F\eta_{x} = (\varepsilon_{Fx})^{-1}$.

We must now check that $\eta$ is natural and that the first
triangle identity holds.  For naturality, we consider the square
below.
\begin{displaymath}
\xymatrix{x \ar[r]^{f} \ar[d]_{\eta_{x}} & y \ar[d]^{\eta_{y}} \\
GFx \ar[r]_{GFf} & GFy }
\end{displaymath}
Applying $F$ to the diagram and using functoriality gives this
square.
\begin{displaymath}
\xymatrix{
Fx \ar[r]^{Ff} \ar[d]_{F\eta_{x}} & Fy \ar[d]^{F\eta_{y}} \\
FGFx \ar[r]_{FGFf} & FGFy }
\end{displaymath}
By the definition of $\eta$, this square is the naturality square
of $(\varepsilon F)^{-1}$ and thus must commute.  By the
faithfulness of $F$, the original square commutes as well and so
$\eta$ is natural.

For the first triangle identity, we consider the composite
\begin{displaymath}
\xymatrix{ Gy \ar[r]^{\eta_{Gy}} & GFGy \ar[r]^{G\varepsilon_{y}}
& Gy. }
\end{displaymath}
Applying $F$ to this yields $FG\varepsilon_{y} \circ F\eta_{Gy}$,
which is by definition $FG \varepsilon_{y} \circ
\varepsilon_{FGy}^{-1}$.  Now the following square commutes by the
naturality of $\varepsilon$.
\begin{displaymath}
\xymatrix{
FGFGy \ar[r]^{\varepsilon_{FGy}} \ar[d]_{FG \varepsilon_{y}} & FGy \ar[d]^{\varepsilon_{y}} \\
FGy \ar[r]_{\varepsilon_{y}} & y }
\end{displaymath}
By the invertibility of $\varepsilon_{y}$, we get that
$\varepsilon_{FGy} = FG\varepsilon_{y}$.  Therefore  $FG
\varepsilon_{y} \circ \varepsilon_{FGy}^{-1} = 1_{FGy}$.  Once
again by the faithfulness of $F$, $G\varepsilon_{y} \circ
\eta_{Gy} = 1_{y}$ and thus the first triangle identity is
satisfied.
\end{proof}

\remark \label{remcat}
We could have just as easily constructed an adjoint equivalence
with $\eta = \beta$ instead of $\varepsilon = \alpha$.  In general
it is not possible to require both of these conditions, though, as
the choice of either $\eta$ or $\epsilon$ fixes the other.
\endremark

It would be possible at this point to prove, by a series of
calculations, an analogous result with $\textbf{Cat}$ replaced by
any bicategory, but instead we choose a different approach that
generalizes more easily to the case of tricategories.

\definition
Let $B$ be a bicategory, and let $f:x \rightarrow y$ be a 1-cell
in $B$.  Then $f$ is an \textit{equivalence} if there exists a
$g:y \rightarrow x$ such that $fg \cong 1_{y}$ and $gf \cong
1_{x}$.
\enddefinition

\lemma\label{functorbicats}
Let $B,C$ be bicategories, and assume that every equivalence $f$
in $C$ is part of an adjoint equivalence $f \dashv_{eq} g$.  Then
every equivalence in the functor bicategory $[B,C]$ is part of an
adjoint equivalence.
\endlemma
\begin{proof}
Let $\alpha:F \Rightarrow G$ be an equivalence in the functor
category.  Then each 1-cell $\alpha_{x}:Fx \rightarrow Gx$ is an
equivalence in $C$, thus we can produce adjoint equivalences
$\alpha_{x} \dashv_{eq} \beta_{x}$ in $C$ for every object $x$ in
$B$.  We will now define a transformation $\beta:F \Rightarrow G$
using the $\beta_{x}$ constructed above as the components on
objects.  Now given a morphism $f:x \rightarrow y$ in $B$, we must
also produce an invertible 2-cell $\beta_{f}: \beta_{y} \circ Gf
\Rightarrow Ff \circ \beta_{x}$ in $C$, subject to the transformation axioms.
We define $\beta_{f}$ by the requirement that it provides the
equality of pasting diagrams given below.
\[
\xy {\ar^{\beta_{x}} (0,0)*+{Gx}; (20,0)*+{Fx} };
{\ar^{\alpha_{x}} (20,0)*+{Fx}; (40,0)*+{Gx} }; {\ar@/^1.5pc/^{1}
(0,0)*+{Gx}; (40,0)*+{Gx} }; {\ar^{Gf} (40,0)*+{Gx};
(40,-13)*+{Gy} }; {\ar_{Gf} (0,0)*+{Gx}; (0,-13)*+{Gy} };
{\ar_{\beta_{y}}  (0,-13)*+{Gy}; (20,-13)*+{Fy} };
{\ar_{\alpha_{y}} (20,-13)*+{Fy}; (40,-13)*+{Gy} }; {\ar^{Ff}
(20,0)*+{Fx}; (20,-13)*+{Fy} }; {\ar@{=>}_{\beta_{f}} (4,-10)*+{};
(16,-3)*+{} }; {\ar@{=>}_{\alpha_{f}} (24,-10)*+{}; (36,-3)*+{} };
{\ar@{=>}_{\varepsilon_{x}} (20,2)*{}; (20,4.3)*{} };
(48,-6.5)*{=}; {\ar_{Gf} (55,0)*+{Gx}; (55,-13)*+{Gy} };
{\ar_{\beta_{y}} (55,-13)*+{Gy}; (75,-13)*+{Fy} };
{\ar_{\alpha_{y}} (75,-13)*+{Fy}; (95,-13)*+{Gy} };
{\ar@/^1.5pc/^{1} (55,-13)*+{Gy}; (95,-13)*+{Gy} }; {\ar^{1}
(55,0)*+{Gx}; (95,0)*+{Gx} }; {\ar^{Gf} (95,0)*+{Gx};
(95,-13)*+{Gy} }; {\ar@{=>}_{\varepsilon_{y}} (75,-10)*{};
(75,-7.5)*{} }; (69,-3)*{\cong}
\endxy
\]
Here, $\varepsilon_{x}$ is the counit of the adjoint equivalence
$\alpha_{x} \dashv_{eq} \beta_{x}$.  This gives a well-defined
$\beta_{f}$ as follows.  Since all the 2-cells in this pasting
diagram are invertible, this equality determines $1_{\alpha_{y}} *
\beta_{f}$.  But $\alpha_{y}$ is an equivalence 1-cell, hence the
functor $\alpha_{y} \circ -$ is an equivalence of categories, so
that  $1_{\alpha_{y}} * \beta_{f}$ determines $\beta_{f}$.

In a similar fashion, we can also construct an invertible 2-cell
$\beta_{f}'$ with the same source and target as $\beta_{f}$ by
requiring it provides the equality of pasting diagrams given
below.
\[
\xy {\ar^{\alpha_{x}} (0,0)*+{Fx}; (20,0)*+{Gx} };
{\ar^{\beta_{x}} (20,0)*+{Gx}; (40,0)*+{Fx} }; {\ar^{Ff}
(40,0)*+{Fx}; (40,-13)*+{Fy} }; {\ar_{Ff} (0,0)*+{Fx};
(0,-13)*+{Fy} }; {\ar_{\alpha_{y}}  (0,-13)*+{Fy}; (20,-13)*+{Gy}
}; {\ar_{\beta_{y}} (20,-13)*+{Gy}; (40,-13)*+{Fy} }; {\ar^{Gf}
(20,0)*+{Gx}; (20,-13)*+{Gy} }; {\ar@/_1.5pc/_{1} (0,-13)*+{Fy};
(40,-13)*+{Fy} }; {\ar@{=>}_{\alpha_{f}} (4,-10)*+{}; (16,-3)*+{}
}; {\ar@{=>}_{\beta_{f}'} (24,-10)*+{}; (36,-3)*+{} };
{\ar@{=>}_{\eta_{y}} (20,-18)*{}; (20,-15.5)*{} }; (48,-6.5)*{=};
{\ar_{Ff} (55,0)*+{Fx}; (55,-13)*+{Fy} }; {\ar^{\alpha_{x}}
(55,0)*+{Fx}; (75,0)*+{Gx} }; {\ar^{\beta_{x}} (75,0)*+{Gx};
(95,0)*+{Fx} }; {\ar@/_1.5pc/_{1} (55,0)*+{Fx}; (95,0)*+{Fx} };
{\ar_{1} (55,-13)*+{Fy}; (95,-13)*+{Fy} }; {\ar^{Ff} (95,0)*+{Fx};
(95,-13)*+{Fy} }; {\ar@{=>}_{\eta_{x}} (75,-5)*{}; (75,-2)*{} };
(69,-10)*{\cong}
\endxy
\]
After applying left and right unit isomorphisms, the pasting
diagram below can be shown to be equal to both $\beta_{f}$ and
$\beta_{f}'$ using the triangle identities, so $\beta_{f} =
\beta_{f}'$.
\[
\xy {\ar^{\beta_{x}} (0,0)*+{Gx}; (20,0)*+{Fx} };
{\ar^{\alpha_{x}} (20,0)*+{Fx}; (40,0)*+{Gx} }; {\ar@/^2pc/^{1}
(0,0)*+{Gx}; (40,0)*+{Gx} }; {\ar^{Gf} (40,0)*+{Gx};
(40,-13)*+{Gy} }; {\ar_{Gf} (0,0)*+{Gx}; (0,-13)*+{Gy} };
{\ar_{\beta_{y}}  (0,-13)*+{Gy}; (20,-13)*+{Fy} };
{\ar_{\alpha_{y}} (20,-13)*+{Fy}; (40,-13)*+{Gy} }; {\ar^{Ff}
(20,0)*+{Fx}; (20,-13)*+{Fy} }; {\ar@{=>}_{\beta_{f}} (4,-10)*+{};
(16,-3)*+{} }; {\ar@{=>}_{\alpha_{f}} (24,-10)*+{}; (36,-3)*+{} };
{\ar@{=>}_{\varepsilon_{x}} (19,3)*{}; (19,6)*{} };
{\ar^{\beta_{x}} (40,0)*+{Gx}; (60,0)*+{Fx} }; {\ar^{Ff}
(60,0)*+{Fx}; (60,-13)*+{Fy} }; {\ar_{\beta_{y}} (40,-13)*+{Gy};
(60,-13)*+{Fy} }; {\ar@{=>}_{\beta_{f}'} (44,-10)*+{}; (56,-3)*+{}
}; {\ar@/_2pc/_{1} (20,-13)*+{Fy}; (60,-13)*+{Fy} };
{\ar@{=>}_{\eta_{y}} (39,-19)*{}; (39,-16)*{} }
\endxy
\]
This defines the components of $\beta$ on both objects and
morphisms.

Now we check that these components satisfy the axioms for a
tranformation.  First, we show that $\beta_{f}$ is natural in $f$.
Given a 2-cell $\delta:f \Rightarrow g$ in $B$, we must show that
\[
\beta_{g} \circ (1_{\beta_{y}} * G\delta) = (F \delta *
1_{\beta_{x}}) \circ \beta_{f}.
\]
This follows from the naturality of both $\alpha_{f}$ and the
coherence isomorphisms used in the definition of $\beta_{f}$.

Second, we must show that $\beta_{1_{x}}$ is the composite
\[
\beta_{x} \circ G1 \cong \beta_{x} 1 \cong \beta_{x} \cong 1
\beta_{x} \cong F1 \circ \beta_{x}
\]
where every isomorphism is given by a unique coherence
isomorphism.  To do this, we need only show that the composite above gives the
equality of pasting diagrams we used to define $\beta_{1_{x}}$.
This is trivial using the unit axiom for the transformation
$\alpha$ and the fact that
\[
\xy {\ar^{\beta_{x}} (0,0)*+{Gx}; (20,0)*+{Fx} };
{\ar^{\alpha_{x}} (20,0)*+{Fx}; (40,0)*+{Gx} }; {\ar^{G1}
(40,0)*+{Gx}; (40,-13)*+{Gx} }; {\ar_{G1} (0,0)*+{Gx};
(0,-13)*+{Gx} }; {\ar_{\beta_{x}} (0,-13)*+{Gx}; (20,-13)*+{Fx} };
{\ar_{\alpha_{x}} (20,-13)*+{Fx}; (40,-13)*+{Gx} }; {\ar^{F1}
(20,0)*+{Fx}; (20,-13)*+{Fx} }; {\ar@{=>}_{\beta_{1_{x}}}
(4,-10)*+{}; (16,-3)*+{} }; {\ar@{=>}_{\alpha_{1_{x}}}
(24,-10)*+{}; (36,-3)*+{} };
\endxy
\]
is the unique coherence isomorphism
\[
(\alpha_{x} \circ \beta_{x}) \circ G1 \cong G1 \circ (\alpha_{x}
\circ \beta_{x})
\]
by the definition of $\beta_{1_{x}}$.

The third and final transformation axiom follows from a similar
proof.
\end{proof}

\lemma
Assume that every equivalence $f$ in $C$ is part of an adjoint
equivalence $f \dashv_{eq} g$, and let $F:B \rightarrow C$ be a
functor which is locally an equivalence.  Then every equivalence
$r$ in $B$ is part of an adjoint equivalence $r \dashv_{eq} s$.
\endlemma
\begin{proof}
If $r:x \rightarrow y$ is an equivalence in $B$, then there is an
$s:y \rightarrow x$ such that $rs \cong 1_{y}$ and $sr \cong
1_{x}$.  Then $Fr$ is an equivalence in $C$ with $Fs$ as a
pseudoinverse.  By hypothesis, we can find an adjoint equivalence
$Fr \dashv_{eq} t$.  By the uniqueness of pseudo-inverses, we must have $t \cong Fs$, so we have an
adjoint equivalence $Fr \dashv_{eq} Fs$.  Since $F$ is locally
full, this means we can find a 2-cell in $B$ which maps to the following composite.
\[
F1 \cong 1 \stackrel{\eta'}{\longrightarrow} FsFr \cong F(sr)
\]
(Here $\eta'$ denotes the unit of the adjoint equivalence $Fr \dashv_{eq} Fs$.)  This 2-cell will be the unit of our adjoint equivalence, and the counit is constructed similarly; the triangle identities follow from coherence for functors and the fact that $F$ is locally faithful.
\end{proof}

\theorem\label{adjeqinbicats}
Let $B$ be a bicategory, and let $f$ be an equivalence in $B$.
Then $f$ is part of an adjoint equivalence $f \dashv_{eq} g$.
\endtheorem
\begin{proof}
Since every equivalence in $\textbf{Cat}$ is part of an adjoint
equivalence, the same is true for $[B^{\textrm{op}},
\textbf{Cat}]$.  Let $Y: B \rightarrow [B^{\textrm{op}},
\textbf{Cat}]$ be the Yoneda embedding.  The functor $Y$ satisfies
the hypotheses of the above lemma, hence there is an adjoint
equivalence $f \dashv_{eq} g$ in $B$.
\end{proof}

\remark
We have actually shown something stronger than the fact that every
equivalence is part of an adjoint equivalence.  We have actually
shown that given any equivalence $f$, a pseudo-inverse $g$, and an
isomorphism $\alpha:fg \cong 1$, there is a unique adjoint
equivalence $f \dashv g$ with $\alpha$ as its counit.  This is
true in $\textbf{Cat}$, hence in any functor bicategory into
$\textbf{Cat}$.  Therefore any two adjoint equivalences $f
\dashv_{eq} g$ with the same counit $\alpha:fg \Rightarrow 1$ will
necessarily have the same unit after applying the Yoneda
embedding, therefore must have the same unit before applying $Y$
since it is a local equivalence.
\endremark

\section{Definitions}
This section will provide the definition of a biadjoint
biequivalence in an arbitrary tricategory $T$.  This proceeds in two steps:  first we define a biadjunction in a tricategory, and then equip it with extra structure to define a biadjoint biequivalence.  There are two
possible options for the definition of a biadjoint biequivalence.
We provide the concise definition first (omitting the ``horizontal
cusp'' axioms), and then explain how it is equivalent to a
definition with additional data and axioms.

Before giving the definition of a biadjunction in a tricategory $T$, we note that our definition is merely the weakening of previous definitions of biadjunctions in \textbf{Gray}-categories.  This weakening is done in the most straightforward manner, and is done because the most natural and concise definition of a biadjunction in a \textbf{Gray}-category uses the \textbf{Gray}-category axioms implicitly.  A different approach to these structures might, for instance, involve constructing the ``free living biadjunction'' - this would be a tricategory $\mathbb{B}$ with the property that biadjunctions in an arbitrary tricategory $T$ would correspond to maps $\mathbb{B} \rightarrow T$.  This is the approach taken by Lack in \cite{lack} in the context of \textbf{Gray}-categories in order to discuss the relationship between biadjunctions and pseudomonads.  Since our focus is on biadjoint biequivalences, and not the more general biadjunctions, we do not proceed in this fashion.

\definition
Let $T$ be a tricategory.  Then a \textit{biadjunction} $f
\dashv_{bi} g$ consists of
\begin{itemize}
\item 1-cells $f:x \rightarrow y, g:y \rightarrow x$, \item
2-cells $\alpha: f \otimes g \Rightarrow I_{y}, \beta:I_{x}
\Rightarrow g \otimes f$, and \item invertible 3-cells $\Phi, \Psi$ below,
\[
\xy {\ar^{r^{\dot}} (0,0)*+{f}; (30,0)*+{f \otimes I} }; {\ar^{1
\otimes \beta} (30,0)*+{f \otimes I}; (60,0)*+{f \otimes (g
\otimes f)} }; {\ar^{a^{\dot}} (60,0)*+{f \otimes (g \otimes f)};
(90,0)*+{(f \otimes g) \otimes f} }; {\ar^{\alpha \otimes 1}
(90,0)*+{(f \otimes g) \otimes f}; (90,-10)*+{I \otimes f} };
{\ar^{l} (90,-10)*+{I \otimes f}; (90,-20)*+{f} }; {\ar_{1}
(0,0)*+{f}; (90,-20)*+{f} }; {\ar@{=>}^{\Phi} (70,-3)*{};
(60,-9)*{} }; {\ar^{l^{\dot}} (0,-35)*+{g}; (30,-35)*+{I \otimes
g} }; {\ar^{\beta \otimes 1} (30,-35)*+{I \otimes g};
(60,-35)*+{(g \otimes f) \otimes g} }; {\ar^{a} (60,-35)*+{(g
\otimes f) \otimes g}; (90,-35)*+{g \otimes (f \otimes g)} };
{\ar^{1 \otimes \alpha} (90,-35)*+{g \otimes (f \otimes g)};
(90,-45)*+{g \otimes I} }; {\ar^{r} (90,-45)*+{g \otimes I};
(90,-55)*+{g} }; {\ar_{1} (0,-35)*+{g}; (90,-55)*+{g} };
{\ar@{=>}^{\Psi} (70,-38)*{}; (60,-44)*{} }
\endxy
\]
\end{itemize}
such that the pasting diagrams in Figures 1 and 2 are both the identity.
\begin{figure}
\[
\xy
{\ar^{\scriptstyle r^{\dot}1} (0,0)*+{\scriptstyle fg};
(15,14.4)*+{\scriptstyle (fI)g} }; {\ar^{\scriptstyle (1 \beta) 1}
(15,14.4)*+{\scriptstyle (fI)g}; (30,28.8)*+{\scriptstyle
(f(gf))g} }; {\ar^{\scriptstyle a^{\dot}1}
(30,28.8)*+{\scriptstyle (f(gf))g}; (45,43.2)*+{\scriptstyle
((fg)f)g} }; {\ar^{\scriptstyle (\alpha 1)1}
(45,43.2)*+{\scriptstyle ((fg)f)g}; (60,57.6)*+{\scriptstyle
(If)g} }; {\ar^{\scriptstyle l1} (60,57.6)*+{\scriptstyle (If)g};
(75,72)*+{\scriptstyle fg} }; {\ar^{\scriptstyle l^{\dot}}
(75,72)*+{\scriptstyle fg}; (90,36)*+{\scriptstyle I(fg)} };
{\ar^{\scriptstyle 1\alpha} (90,36)*+{\scriptstyle I(fg)};
(105,0)*+{\scriptstyle II} }; {\ar_{\scriptstyle 1 l^{\dot}}
(0,0)*+{\scriptstyle fg}; (15,-14.4)*+{\scriptstyle f(Ig)} };
{\ar_{\scriptstyle 1(\beta 1)} (15,-14.4)*+{\scriptstyle f(Ig)};
(30,-28.8)*+{\scriptstyle f((gf)g)} }; {\ar_{\scriptstyle 1a}
(30,-28.8)*+{\scriptstyle f((gf)g)}; (45,-43.2)*+{\scriptstyle
f(g(fg))} }; {\ar_{\scriptstyle 1(1 \alpha)}
(45,-43.2)*+{\scriptstyle f(g(fg))}; (60,-57.6)*+{\scriptstyle
f(gI)} }; {\ar_{\scriptstyle 1r} (60,-57.6)*+{\scriptstyle f(gI)};
(75,-72)*+{\scriptstyle fg} }; {\ar_{\scriptstyle r^{\dot}}
(75,-72)*+{\scriptstyle fg}; (90,-36)*+{\scriptstyle (fg)I} };
{\ar_{\scriptstyle \alpha 1} (90,-36)*+{\scriptstyle (fg)I};
(105,0)*+{\scriptstyle II} }; {\ar@/^2pc/^{\scriptstyle l}
(105,0)*+{\scriptstyle II}; (130,0)*+{\scriptstyle I} };
{\ar@/_2pc/_{\scriptstyle r} (105,0)*+{\scriptstyle II};
(130,0)*+{\scriptstyle I} }; {\ar^{\scriptstyle a}
(15,14.4)*+{\scriptstyle (fI)g}; (15,-14.4)*+{\scriptstyle f(Ig)}
}; {\ar@/^1pc/^{\scriptstyle a} (30,28.8)*+{\scriptstyle
(f(gf))g}; (30,-28.8)*+{\scriptstyle f((gf)g)} };
{\ar_{\scriptstyle a} (45,43.2)*+{\scriptstyle ((fg)f)g};
(65,0)*+{\scriptstyle (fg)(fg)} }; {\ar^{\scriptstyle a^{\dot}}
(45,-43.2)*+{\scriptstyle f(g(fg))}; (65,0)*+{\scriptstyle
(fg)(fg)} }; {\ar^{\scriptstyle a^{\dot}}
(60,-57.6)*+{\scriptstyle f(gI)}; (90,-36)*+{\scriptstyle (fg)I}
}; {\ar_{\scriptstyle a} (60,57.6)*+{\scriptstyle (If)g};
(90,36)*+{\scriptstyle I(fg)} }; {\ar_{\scriptstyle \alpha 1}
(65,0)*+{\scriptstyle (fg)(fg)}; (90,36)*+{\scriptstyle I(fg)} };
{\ar^{\scriptstyle 1 \alpha} (65,0)*+{\scriptstyle (fg)(fg)};
(90,-36)*+{\scriptstyle (fg)I} };
{\ar@/^5pc/^{\scriptstyle 1} (0,0)*+{\scriptstyle fg}; (75,72)*+{\scriptstyle fg} };
{\ar@/_5pc/_{\scriptstyle 1} (0,0)*+{\scriptstyle fg};
(75,-72)*+{\scriptstyle fg} }; {\ar@/^2pc/^{\scriptstyle \alpha}
(75,72)*+{\scriptstyle fg}; (130,0)*+{\scriptstyle I} };
{\ar@/_2pc/_{\scriptstyle \alpha} (75,-72)*+{\scriptstyle fg};
(130,0)*+{\scriptstyle I} }; (9,0)*{\scriptstyle \Downarrow \mu};
(24,0)*{\scriptstyle \cong}; (45,0)*{\scriptstyle \Downarrow \pi};
(65,28.8)*{\scriptstyle \cong}; (65,-28.8)*{\scriptstyle \cong};
(73,57.6)*{\scriptstyle \Downarrow \lambda};
(73,-57.6)*{\scriptstyle \Downarrow \rho}; (88,0)*{\scriptstyle
\cong}; (117,0)*{\scriptstyle \cong}; (105,30.4)*{\scriptstyle
\cong}; (105,-30.4)*{\scriptstyle \cong};  (28,-44)*{\scriptstyle
\Downarrow 1 \Psi}; (28,44)*{\scriptstyle \Downarrow \Phi^{-1} 1}
\endxy
\]
\caption{First pasting}
\end{figure}
\begin{figure}
\[
\xy {\ar@/^2pc/^{\scriptstyle \beta} (0,0)*+{\scriptstyle I};
(55,72)*+{\scriptstyle gf} }; {\ar@/_2pc/_{\scriptstyle \beta}
(0,0)*+{\scriptstyle I}; (55,-72)*+{\scriptstyle gf} };
{\ar@/^2pc/^{\scriptstyle l^{\dot}} (0,0)*+{\scriptstyle I};
(25,0)*+{\scriptstyle II} }; {\ar@/_2pc/_{\scriptstyle r^{\dot}}
(0,0)*+{\scriptstyle I}; (25,0)*+{\scriptstyle II} };
{\ar^{\scriptstyle 1 \beta} (25,0)*+{\scriptstyle II};
(40,36)*+{\scriptstyle I(gf)} }; {\ar_{\scriptstyle \beta1}
(25,0)*+{\scriptstyle II}; (40,-36)*+{\scriptstyle (gf)I} };
{\ar^{\scriptstyle l} (40,36)*+{\scriptstyle I(gf)};
(55,72)*+{\scriptstyle gf} }; {\ar_{\scriptstyle r}
(40,-36)*+{\scriptstyle (gf)I}; (55,-72)*+{\scriptstyle gf} };
{\ar_{\scriptstyle \beta1} (40,36)*+{\scriptstyle I(gf)};
(65,0)*+{\scriptstyle (gf)(gf)} }; {\ar^{\scriptstyle 1 \beta}
(40,-36)*+{\scriptstyle (gf)I}; (65,0)*+{\scriptstyle (gf)(gf)} };
{\ar^{\scriptstyle l^{\dot}1} (55,72)*+{\scriptstyle gf};
(70,57.6)*+{\scriptstyle (Ig)f} }; {\ar^{\scriptstyle
(\beta^{\dot}1)1} (70,57.6)*+{\scriptstyle (Ig)f};
(85,43.2)*+{\scriptstyle ((gf)g)f} }; {\ar^{\scriptstyle a1}
(85,43.2)*+{\scriptstyle ((gf)g)f}; (100,28.8)*+{\scriptstyle
(g(fg))f} }; {\ar^{\scriptstyle (1 \alpha)1}
(100,28.8)*+{\scriptstyle (g(fg))f}; (115,14.4)*+{\scriptstyle
(gI)f} }; {\ar^{\scriptstyle r1} (115,14.4)*+{\scriptstyle (gI)f};
(130,0)*+{\scriptstyle gf} }; {\ar_{\scriptstyle 1 r^{\dot}}
(55,-72)*+{\scriptstyle gf}; (70, -57.6)*+{\scriptstyle g(fI)} };
{\ar_{\scriptstyle 1(1 \beta)} (70,-57.6)*+{\scriptstyle g(fI)};
(85,-43.2)*+{\scriptstyle g(f(gf))} }; {\ar_{\scriptstyle 1
a^{\dot}} (85,-43.2)*+{\scriptstyle g(f(gf))};
(100,-28.8)*+{\scriptstyle g((fg)f)} }; {\ar_{\scriptstyle
1(\alpha 1)} (100,-28.8)*+{\scriptstyle g((fg)f)};
(115,-14.4)*+{\scriptstyle g(If)} }; {\ar_{\scriptstyle 1l}
(115,-14.4)*+{\scriptstyle g(If)}; (130,0)*+{\scriptstyle gf} };
{\ar_{\scriptstyle a} (115,14.4)*+{\scriptstyle (gI)f};
(115,-14.4)*+{\scriptstyle g(If)} }; {\ar@/_1pc/_{\scriptstyle a}
(100,28.8)*+{\scriptstyle (g(fg))f}; (100,-28.8)*+{\scriptstyle
g((fg)f)} }; {\ar_{\scriptstyle a^{\dot}} (65,0)*+{\scriptstyle
(gf)(gf)}; (85,43.2)*+{\scriptstyle ((gf)g)f} };
{\ar^{\scriptstyle a} (65,0)*+{\scriptstyle (gf)(gf)};
(85,-43.2)*+{\scriptstyle g(f(gf))} }; {\ar_{\scriptstyle
a^{\dot}} (40,36)*+{\scriptstyle I(gf)}; (70,57.6)*+{\scriptstyle
(Ig)f} }; {\ar^{\scriptstyle a} (40,-36)*+{\scriptstyle (gf)I};
(70, -57.6)*+{\scriptstyle g(fI)} };
{\ar@/^5pc/^{\scriptstyle 1} (55,72)*+{\scriptstyle gf}; (130,0)*+{\scriptstyle gf} };
{\ar@/_5pc/_{\scriptstyle 1} (55,-72)*+{\scriptstyle gf};
(130,0)*+{\scriptstyle gf} }; (12,0)*{\scriptstyle \cong};
(121,0)*{\scriptstyle \Downarrow \mu}; (106,0)*{\scriptstyle
\cong}; (85,0)*{\scriptstyle \Downarrow \pi};
(65,28.8)*{\scriptstyle \cong}; (65,-28.8)*{\scriptstyle \cong};
(57,57.6)*{\scriptstyle \Downarrow \lambda};
(57,-57.6)*{\scriptstyle \Downarrow \rho}; (42,0)*{\scriptstyle
\cong}; (25,30.4)*{\scriptstyle \cong}; (25,-30.4)*{\scriptstyle
\cong};  (102,-44)*{\scriptstyle \Downarrow 1 \Phi}; (102, 44)*{\scriptstyle \Downarrow \Psi^{-1} 1}
\endxy
\]
\caption{Second pasting} \end{figure}

\enddefinition

\remark
In the presence of the simplifying assumption that the tricategory
$T$ is actually a strict, cubical tricategory (i.e., a
$\mathbf{Gray}$-category), the axioms simplify to the
equality of pasting diagrams below.
\[
\xy {\ar^{f \beta g} (0,0)*+{fg}; (10,14)*+{fgfg} }; {\ar^{\alpha
fg} (10,14)*+{fgfg}; (30,14)*+{fg} }; {\ar^{\alpha} (30,14)*+{fg};
(40,0)*+{I} }; {\ar_{1} (0,0)*+{fg}; (20,0)*+{fg} }; {\ar_{\alpha}
(20,0)*+{fg}; (40,0)*+{I} }; {\ar^{fg \alpha} (10,14)*+{fgfg};
(20,0)*+{fg} }; (10,5)*{\scriptstyle \Downarrow f \Psi};
(27,7)*{\cong}; {\ar^{f \beta g} (60,0)*+{fg}; (70,14)*+{fgfg} };
{\ar^{\alpha fg} (70,14)*+{fgfg}; (90,14)*+{fg} }; {\ar^{\alpha}
(90,14)*+{fg}; (100,0)*+{I} }; {\ar_{1} (60,0)*+{fg}; (80,0)*+{fg}
}; {\ar_{\alpha} (80,0)*+{fg}; (100,0)*+{I} }; {\ar_{1}
(60,0)*+{fg}; (90,14)*+{fg} }; (72,9)*{\scriptstyle \Downarrow
\Phi g}; (90,5)*{=}; {\ar@{=} (45,7)*{}; (55,7)*{} }
\endxy
\]
\[
\xy {\ar^{ \beta } (0,0)*+{I}; (10,14)*+{gf} }; {\ar^{\beta gf}
(10,14)*+{gf}; (30,14)*+{gfgf} }; {\ar^{g \alpha f}
(30,14)*+{gfgf}; (40,0)*+{gf} }; {\ar_{\beta} (0,0)*+{I};
(20,0)*+{gf} }; {\ar_{1} (20,0)*+{gf}; (40,0)*+{gf} }; {\ar^{gf
\beta} (20,0)*+{gf}; (30,14)*+{gfgf} }; (13,7)*{\cong};
(30,5)*{\scriptstyle \Downarrow g \Phi}; {\ar^{ \beta }
(60,0)*+{I}; (70,14)*+{gf} }; {\ar^{\beta gf} (70,14)*+{gf};
(90,14)*+{gfgf} }; {\ar^{g \alpha f} (90,14)*+{gfgf};
(100,0)*+{gf} }; {\ar_{\beta} (60,0)*+{I}; (80,0)*+{gf} };
{\ar_{1} (80,0)*+{gf}; (100,0)*+{gf} }; {\ar_{1} (70,14)*+{gf};
(100,0)*+{gf} }; (88,9)*{\scriptstyle \Downarrow \Psi f};
(70,5)*{=}; {\ar@{=} (45,5)*{}; (55,5)*{} }
\endxy
\]
See \cite{st}, \cite{ver}, or \cite{lack} for earlier definitions.
\endremark

\definition
Let $T$ be a tricategory.  Then a \textit{biadjoint biequivalence}
$f \dashv_{bieq} g$ consists of
\begin{itemize}
\item a biadjunction $f \dashv_{bi} g$  and \item adjoint equivalences $\alpha
\dashv_{eq} \alpha^{\dot}, \beta \dashv_{eq} \beta^{\dot}$ in the
respective hom-bicategories.
\end{itemize}
\enddefinition

It is also possible to give a longer version of a biadjoint
biequivalence which includes extra data satisfying the so-called
horizontal cusp axioms.  Such a definition is equivalent to the
one given above by using the calculus of mates, as we explain
below.

The extra data needed to express the horizontal cusp axioms are a
pair of invertible 3-cells $\underline{\Phi}, \underline{\Psi}$.
\[
\xy {\ar^{l^{\dot}} (0,0)*+{f}; (30,0)*+{I \otimes f} };
{\ar^{\alpha^{\dot} \otimes 1} (30,0)*+{I \otimes f}; (60,0)*+{(f
\otimes g) \otimes f} }; {\ar^{a} (60,0)*+{(f \otimes g) \otimes
f}; (90,0)*+{f \otimes (g \otimes f)} }; {\ar^{1 \otimes
\beta^{\dot}} (90,0)*+{f \otimes (g \otimes f)}; (90,-10)*+{f
\otimes I} }; {\ar^{r} (90,-10)*+{f \otimes I}; (90,-20)*+{f} };
{\ar_{1} (0,0)*+{f}; (90,-20)*+{f} }; {\ar@{=>}^{\underline{\Phi}}
(70,-3)*{}; (60,-9)*{} }; {\ar^{r^{\dot}} (0,-35)*+{g};
(30,-35)*+{g \otimes I} }; {\ar^{1 \otimes \alpha^{\dot}}
(30,-35)*+{g \otimes I}; (60,-35)*+{g \otimes (f \otimes g)} };
{\ar^{a^{\dot}} (60,-35)*+{g \otimes (f \otimes g)}; (90,-35)*+{(g
\otimes f) \otimes g} }; {\ar^{\beta^{\dot} \otimes 1}
(90,-35)*+{(g \otimes f) \otimes g}; (90,-45)*+{I \otimes g} };
{\ar^{l} (90,-45)*+{I \otimes g}; (90,-55)*+{g} }; {\ar_{1}
(0,-35)*+{g}; (90,-55)*+{g} }; {\ar@{=>}^{\underline{\Psi}}
(70,-38)*{}; (60,-44)*{} }
\endxy
\]

These additional 3-cells are then required to satisfy the
horizontal cusp axioms, named for their relationship with certain
``braid movie moves'' between braided surfaces in $\mathbb{R}^{4}$.
One such axiom, written categorically, is given below.

\[
\xy {\ar^{1 \beta^{\dot}} (0,0)*+{f(gf)}; (15,0)*+{f1} }; {\ar^{1
\beta} (15,0)*+{f1}; (25,-10)*+{f(gf)} }; {\ar@/_0.5pc/_{1}
(0,0)*+{f(gf)}; (25,-10)*+{f(gf)} }; {\ar_{a^{\dot}}
(25,-10)*+{f(gf)}; (40,-10)*+{(fg)f} }; {\ar_{\alpha 1}
(40,-10)*+{(fg)f}; (55,-10)*+{1f} }; {\ar_{l} (55,-10)*+{1f};
(70,0)*+{f} }; {\ar^{r} (15,0)*+{f1}; (70,0)*+{f} };
(12,-4)*{\scriptstyle \Downarrow \Delta}; (42.5,-5)*{\Downarrow
\Phi^{-1}}; {\ar_{a^{\dot}} (0,-30)*+{f(gf)}; (10,-40)*+{(fg)f} };
{\ar_{\alpha 1} (10,-40)*+{(fg)f}; (20,-50)*+{1f} }; {\ar^{1}
(10,-40)*+{(fg)f}; (30,-40)*+{(fg)f} }; {\ar_{\alpha^{\dot} 1}
(20,-50)*+{1f}; (30,-40)*+{(fg)f} }; {\ar^{a} (30,-40)*+{(fg)f};
(45,-40)*+{f(gf)} }; {\ar^{1 \beta^{\dot}} (45,-40)*+{f(gf)};
(60,-40)*+{f1} }; {\ar^{r}  (60,-40)*+{f1}; (70,-50)*+{f} };
{\ar_{l} (20,-50)*+{1f}; (70,-50)*+{f} }; {\ar@/^0.5pc/^{1}
(0,-30)*+{f(gf)}; (45,-40)*+{f(gf)} }; (20,-35)*{\Downarrow
\eta_{a}}; (20,-44)*{\scriptstyle \Downarrow \underline{\Gamma}
1}; (45,-45)*{\Downarrow \underline{\Phi}}; {\ar@{=} (35,-18)*{};
(35,-28)*{} }
\endxy
\]

It is then clear that this axiom merely says that
$\underline{\Phi}$ is the mate of $\Phi$, and similarly for
$\underline{\Psi}$ and $\Psi$.  All of the horizontal cusp-type
axioms can be expressed in this fashion.



\section{Biequivalences in \textbf{Bicat}}
This section presents a computational proof that every
biequivalence in the tricategory $\textbf{Bicat}$ is part of a
biadjoint biequivalence.  The proof here will proceed much as the
proof in $\textbf{Cat}$ did, by relying on an alternate
description of biequivalences.  Thus we begin with a simple lemma.

\lemma
Let $F: B \rightarrow C$ be a functor between bicategories.  Then
$F$ is biessentially surjective and locally an equivalence of
categories if and only if there is a functor $G:C \rightarrow B$
such that $FG \simeq 1_{C}$ and $GF \simeq 1_{B}$ in the
respective functor bicategories.
\endlemma
\begin{proof}
Since $F$ is biessentially surjective, for every object $c$ in $C$
we can find an object $b$ in $B$ and an adjoint equivalence $f_{c}
\dashv_{eq} g_{c}$ between $Fb$ and $c$ by Theorem \ref{adjeqinbicats}; here we
choose $f_{c}$ to have source $Fb$ and target $c$.  We choose such
an adjoint equivalence for every $c$, and define the functor $G$
on objects by $Gc = b$. Now $F_{b,b'}:B(b,b') \rightarrow C(Fb,
Fb')$ is an equivalence of categories for every pair $b,b'$, and
we choose an adjoint equivalence $F_{b,b'} \dashv_{eq} \tilde{G}_{b,b'}$.
We define the functor $G$ on hom-categories $G_{b,b'}:C(c, c')
\rightarrow B(Gc,Gc')$ to be the composite
\[
C(c,c') \stackrel{f_{c}^{*}}{\rightarrow} C(Fb, c')
\stackrel{g_{c*}}{\rightarrow} C(Fb, Fb')
\stackrel{\tilde{G}_{b,b'}}{\rightarrow} B(b,b').
\]

We must now construct isomorphisms $1_{Gc} \cong G(1_{c})$, $Gf
\circ Gg \cong G(f \circ g)$ and check the axioms for a functor.
For the first of these, we compute that
\[
G(1_{c}) = \tilde{G}_{b,b} \Big( g_{c} \circ (1_{c} \circ f_{c}) \Big),
\]
while $1_{Gc} = 1_{b}$.  Now note that the adjoint equivalence
$F_{b,b'} \dashv_{eq} \tilde{G}_{b,b'}$ has a unit isomorphism $1
\Rightarrow \tilde{G}_{b,b'} \circ F_{b,b'}$, and when specialized to the
case $b=b'$ and then evaluated at $1_{b}$ yields
\[
1_{b} \cong \tilde{G}_{b,b}( F_{b,b}(1_{b})).
\]
Since $F$ is a functor, we have an isomorphism $\varphi_{0}^{F}:
1_{Fb} \cong F(1_{b})$ which we can compose with the previous
isomorphism to get
\[
1_{b} \cong \tilde{G}_{b,b}( F_{b,b}(1_{b}))
\stackrel{\tilde{G}_{b,b}(\varphi_{0}^{F})}{\longrightarrow} \tilde{G}_{b,b}(1_{Fb})
\]
which we denote by $\overline{\varphi}_{0}^{G}$.  Writing
$\eta_{c}: 1_{Fb} \Rightarrow g_{c} \circ f_{c}$ for the unit of
the adjoint equivalence $f_{c} \dashv_{eq} g_{c}$, we obtain the
isomorphism
\[
\varphi_{0}^{G}: 1_{b} \cong G(1_{c})
\]
as the following composite.
\[
1_{b} \stackrel{\overline{\varphi}_{0}^{G}}{\longrightarrow} \tilde{G}_{b,b}(1_{Fb})
\stackrel{\tilde{G}_{b,b}(\eta_{c})}{\longrightarrow} \tilde{G}_{b,b}(g_{c} \circ f_{c})
\stackrel{\tilde{G}_{b,b}(1*l^{-1})}{\longrightarrow} \tilde{G}_{b,b} \Big( g_{c} \circ (1
\circ f_{c}) \Big)
\]

The isomorphism $\varphi_{2}^{G}: Gf \circ Gg \cong G(f \circ g)$
is obtained in a similar fashion, and the functor axioms for $G$
follow from those for $F$ and the adjoint equivalence axioms for
both $f_{c} \dashv_{eq} g_{c}$ and $F_{b,b'} \dashv_{eq}
\tilde{G}_{b,b'}$.

We must finally check that $FG \simeq 1_{C}$ and $GF \simeq 1_{B}$
in the relevant functor bicategories.  For the first of these,
note that $FG(c) = Fb$ by construction. We already have equivalence
1-cells $f_{c}:FG(c) \rightarrow 1_{C}(c)$ that we now need to
complete to a natural transformation.  This requires that we give
natural isomorphisms
\[
f_{r}: f_{c'}g_{c'}rf_{c} \cong rf_{c},
\]
one for each $r$, which we define to be the obvious whiskering of the counit
isomorphism $f_{c'}g_{c'} \cong 1_{c'}$ composed with the left
unit isomorphism.  The adjoint equivalence axioms for $f_{c}
\dashv_{eq} g_{c}$ and coherence for bicategories imply all of the
transformation axioms.  This shows that $FG \simeq 1_{C}$, and we
leave it to the reader to prove $GF \simeq 1_{B}$.
\end{proof}

\theorem\label{biadjbieq1}
Let $F:B \rightarrow C$ be a biequivalence between bicategories.
Then there is a biadjoint biequivalence $F \dashv_{bieq} G$.
\endtheorem
\begin{proof}
Since $F$ is a biequivalence, choose a functor $G:C \rightarrow B$
such that $FG$ is equivalent to $1_{C}$ in the bicategory
$\textbf{Bicat}(C,C)$ and $GF$ is equivalent to $1_{B}$ in
$\textbf{Bicat}(B,B)$. Taking any equivalence $\alpha:FG
\Rightarrow 1_{C}$ exhibiting this fact, we can construct an
adjoint equivalence $\alpha \dashv_{eq} \alpha^{\dot}$ in
$\textbf{Bicat}(C,C)$ by Theorem \ref{adjeqinbicats}. We will write $\Gamma:\alpha
\alpha^{\dot} \Rrightarrow 1$ and $\underline{\Gamma}:1
\Rrightarrow \alpha^{\dot} \alpha$ for the counit and unit of this
adjoint equivalence, respectively.

Now we construct the adjoint equivalence $\beta \dashv_{eq}
\beta^{\dot}$ (between $1_{B}$ and $GF$) and the invertible
modification $\Phi$ simultaneously.  The component of $\beta$ at
an object $b \in B$ is a 1-cell $\beta_{b}:b \rightarrow GFb$. The
component of $\Phi$ at $b \in B$ is an invertible 2-cell in $C$
\[
\xy {\ar^{1} (0,0)*+{Fb}; (20,0)*+{Fb} }; {\ar^{F \beta_{b}}
(20,0)*+{Fb}; (40,0)*+{FGFb} }; {\ar^{1} (40,0)*+{FGFb};
(60,0)*+{FGFb} }; {\ar^{\alpha_{Fb}} (60,0)*+{FGFb};
(60,-10)*+{Fb} }; {\ar^{1} (60,-10)*+{Fb}; (60,-20)*+{Fb} };
{\ar_{1} (0,0)*+{Fb}; (60,-20)*+{Fb} }; {\ar@{=>}_{\Phi}
(45,-4)*{}; (35,-10)*{} }
\endxy
\]
since the associativity and unit 2-cells in $\textbf{Bicat}$ have
identities as their components.  By coherence, such an invertible
2-cells determines and is determined by an invertible 2-cell
$\widetilde{\Phi}:\alpha_{Fb} \circ F(\beta_{b}) \Rightarrow
1_{Fb}$. Since $\alpha$ is an equivalence, giving such an
isomorphism is equivalent to giving an isomorphism
$\alpha^{\dot}_{Fb} \cong F(\beta_{b})$. Now $F$ is locally an
equivalence of categories, and in particular essentially
surjective, so there exists a morphism $\beta_{b}:b \rightarrow
GFb$ such that $F(\beta_{b}) \cong \alpha_{Fb}^{\dot}$.  For every
object $b \in B$, choose such a $\beta_{b}$ and a specified
isomorphism $\delta_{b}:F(\beta_{b}) \cong \alpha_{Fb}^{\dot}$.

For the component of $\beta$ at $f:b \rightarrow c$, consider the
following composite.
\[
\begin{array}{rcccl}
F(\beta_{c} \circ f) & \cong & F\beta_{c} \circ Ff & \stackrel{\delta_{c}*1}{\rightarrow} & \alpha_{Fc}^{\dot} \circ Ff \\
& \stackrel{\alpha_{f}^{\dot}}{\rightarrow} & FGFf \circ \alpha_{Fb}^{\dot} & \stackrel{1 * \delta_{b}^{-1}}{\rightarrow} & FGFf \circ F\beta_{b} \\
& & & \cong & F \Big( GFf \circ \beta_{b} \Big)
\end{array}
\]
Since $F$ is locally an equivalence, there is a unique isomorphism
\[
\beta_{f}:\beta_{c} \circ f \Rightarrow GFf \circ \beta_{b}
\]
that maps to the composite above.  It is then simple to check that
$\beta$ is a transformation $1 \Rightarrow GF$, and that it is an
equivalence.  This construction also immediately implies that
$\delta$ is an invertible modification
\[
\delta: 1_{F} \otimes \beta \Rrightarrow \alpha^{\dot}.
\]
We then define the adjoint equivalence $\beta \dashv_{eq}
\beta^{\dot}$ to be any adjoint equivalence containing $\beta$.

The 2-cell $\widetilde{\Phi}: \alpha_{Fb} \circ F(\beta_{b})
\Rightarrow 1_{Fb}$ is defined to be the following composite.
\[
\alpha_{Fb} \circ F(\beta_{b}) \stackrel{1 *
\delta_{b}}{\Rightarrow} \alpha_{Fb} \circ \alpha_{Fb}^{\dot}
\stackrel{\Gamma}{\Rightarrow} 1_{Fb}
\]
By coherence for bicategories, this determines the 2-cell
$\Phi_{b}$ uniquely.  These 2-cells $\Phi_{b}$ then give the data
for an invertible modification $\Phi$ since all of the cells used
to construct the $\Phi_{b}$ are either components of modifications
or are appropriately natural.

All that remains is to construct the invertible 3-cell $\Psi$ and
to check the two biadjunction axioms.  Before doing so, we remind the reader
that, for 3-cells in a tricategory, $\circ$ denotes the composition
along 2-cell boundaries, $*$ denotes composition along 1-cell
boundaries, and $\otimes$ denotes composition along 0-cell
boundaries.  Now the second axiom determines the 3-cell
$(1_{G} \otimes \Psi) * 1_{\beta}$.  Since the 2-cell $\beta$ is an equivalence, the
functor $- \circ \beta$ is an equivalence of categories, and in
particular the cell $(1_{G} \otimes \Psi) * 1_{\beta}$ uniquely
determines the cell $1_{G} \otimes \Psi$.  Similarly, since the
functor $F$ is a biequivalence, $G$ is also, so the functor $G
\circ -$ is a biequivalence of bicategories; thus $1_{G} \otimes
\Psi$ uniquely determines the invertible modification $\Psi$.  By
construction, the second biadjoint biequivalence axiom is
satisfied.

Now we show that this choice of $\Psi$ satisfies the first
biadjoint biequivalence axiom.  First, note that, while
$\textbf{Bicat}$ is not a $\textbf{Gray}$-category, it does have a
strictly associative and unital composition law for 1-cells in the
following sense.  The composite $H(GF)$ equals the composite
$(HG)F$, and similarly $F1 = F = 1F$, but we still have
associativity and unit equivalences for this composition law.
These are the 1-cells in the biadjoint biequivalence axioms
labeled $a, l, r$, and they have components on objects given by
identites and components on morphisms given by unique coherence
2-cells.  For examples, the transformation $l: 1F \Rightarrow F$
has its component at an object $x$ the identity $1_{Fx}:Fx
\rightarrow Fx$ and its component at a 1-cell $f:x \rightarrow y$
the unique coherence cell
\[
1_{Fy} \circ Ff \cong Ff \circ 1_{Fx}.
\]
Similarly, the modifications $\pi, \mu, \lambda,$ and $\rho$ in
$\textbf{Bicat}$ all have unique coherence 2-cells as their
components.  Thus coherence for bicategories reduces the first
biadjoint biequivalence to checking that the pasting
\[
\xy {\ar^{1} (0,0)*+{FGx}; (30,15)*+{FGx} }; {\ar^{\alpha_{x}}
(30,15)*+{FGx}; (60,0)*+{x} }; {\ar_{1} (0,0)*+{FGx};
(30,-15)*+{FGx} }; {\ar_{\alpha_{x}} (30,-15)*+{FGx}; (60,0)*+{x}
}; {\ar_{F\beta_{Gx}} (0,0)*+{FGx}; (30,0)*+{FGFGx} };
{\ar^{FG\alpha_{x}} (30,0)*+{FGFGx}; (30,-15)*+{FGx} };
{\ar_{\alpha_{FGx}} (30,0)*+{FGFGx}; (30,15)*+{FGx} };
(23,6)*{\scriptstyle \Downarrow \tilde{\Phi}^{-1}_{Gx}}; (23,-6)*{\scriptstyle \Downarrow
F\tilde{\Psi}_{x}}; (45,0)*{\cong_{\alpha}}
\endxy
\]
is equal to the identity on $\alpha \circ 1_{FG}$.  (Here we use
the same convention that $\widetilde{\Psi}$ is derived from $\Psi$
via unique coherence isomorphisms.)  From this point on, we mark
our naturality isomorphisms with a subscript to indicate which
transformation they are naturality isomorphisms for to avoid
confusion, and we refer to all instances of the above pasting
diagram as ``Axiom 1'', perhaps with some descriptor to indicate
which object $x$ is being used.

First, note that Axiom 1 is the identity if and only if it is
the identity when $x$ is of the form $Fy$ for some $y$ in $B$.
Indeed, consider the following pasting diagram.
\[
\xy {\ar^{1} (0,0)*+{FGx}; (30,15)*+{FGx} }; {\ar^{\alpha_{x}}
(30,15)*+{FGx}; (60,0)*+{x} }; {\ar_{1} (0,0)*+{FGx};
(30,-15)*+{FGx} }; {\ar_{\alpha_{x}} (30,-15)*+{FGx}; (60,0)*+{x}
}; {\ar_{F\beta_{Gx}} (0,0)*+{FGx}; (30,0)*+{FGFGx} };
{\ar^{FG\alpha_{x}} (30,0)*+{FGFGx}; (30,-15)*+{FGx} };
{\ar_{\alpha_{FGx}} (30,0)*+{FGFGx}; (30,15)*+{FGx} };
(23,6)*{\scriptstyle \Downarrow \tilde{\Phi}^{-1}_{Gx}};
(23,-6)*{\scriptstyle \Downarrow F\tilde{\Psi}_{x}};
(45,0)*{\cong_{\alpha}}; {\ar_{FG \alpha_{x}} (0,15)*+{FGFGx};
(0,0)*+{FGx} }; {\ar^{1} (0,15)*+{FGFGx}; (30,30)*+{FGFGx} };
{\ar^{\alpha_{FGx}} (30,30)*+{FGFGx}; (60,15)*+{FGx} };
{\ar^{\alpha_{x}} (60,15)*+{FGx}; (60,0)*+{x} }; {\ar^{FG
\alpha_{x}} (30,30)*+{FGFGx}; (30,15)*+{FGx} };
(16,18)*{\cong_{1}}; (45,18)*{\cong_{\alpha}}
\endxy
\]
Using the modification and transformation axioms, it is equal to
the pasting below.
\[
\xy {\ar^{1} (0,0)*+{FGFGx}; (30,15)*+{FGFGx} };
{\ar^{\alpha_{FGx}} (30,15)*+{FGFGx}; (60,0)*+{FGx} }; {\ar_{1}
(0,0)*+{FGFGx}; (30,-15)*+{FGFGx} }; {\ar_{\alpha_{FGx}}
(30,-15)*+{FGFGx}; (60,0)*+{FGx} }; {\ar_{F\beta_{GFGx}}
(0,0)*+{FGFGx}; (30,0)*+{FGFGFGx} }; {\ar^{FG\alpha_{FGx}}
(30,0)*+{FGFGFGx}; (30,-15)*+{FGFGx} }; {\ar_{\alpha_{FGFGx}}
(30,0)*+{FGFGFGx}; (30,15)*+{FGFGx} }; (23,6)*{\scriptstyle
\Downarrow \tilde{\Phi}^{-1}_{GFGx}}; (23,-6)*{\scriptstyle
\Downarrow F\tilde{\Psi}_{FGx}}; (45,0)*{\cong_{\alpha}}; {\ar_{FG
\alpha_{x}} (0,0)*+{FGFGx}; (0,-15)*+{FGx} }; {\ar_{1}
(0,-15)*+{FGx}; (30,-30)*+{FGx} }; {\ar^{FG\alpha_{x}}
(30,-15)*+{FGFGx}; (30,-30)*+{FGx} }; {\ar^{\alpha_{x}}
(60,0)*+{FGx}; (60,-15)*+{x} }; {\ar_{\alpha_{x}} (30,-30)*+{FGx};
(60,-15)*+{x} }; (16,-16)*{\cong_{1}}; (45,-16)*{\cong_{\alpha}}
\endxy
\]
Thus the pasting Axiom 1 for $x$ is the identity if and only if
Axiom 1 for $FGx$ is the identity, so taking $y=Gx$ proves the
claim.  From this point, we replace $x$ with $Fy$.

Now Axiom 1 for $Fy$ is the identity if and only if the following
pasting diagram is the identity since $F\beta_{y}$ is an
equivalence 1-cell.
\[
\xy {\ar^{1} (0,0)*+{FGFy}; (30,15)*+{FGFy} }; {\ar^{\alpha_{Fy}}
(30,15)*+{FGFy}; (60,0)*+{Fy} }; {\ar_{1} (0,0)*+{FGFy};
(30,-15)*+{FGFy} }; {\ar_{\alpha_{Fy}} (30,-15)*+{FGFy};
(60,0)*+{Fy} }; {\ar_{F\beta_{GFy}} (0,0)*+{FGFy};
(30,0)*+{FGFGFy} }; {\ar^{FG\alpha_{Fy}} (30,0)*+{FGFGFy};
(30,-15)*+{FGFy} }; {\ar_{\alpha_{FGFy}} (30,0)*+{FGFGFy};
(30,15)*+{FGFy} }; (23,6)*{\scriptstyle \Downarrow
\tilde{\Phi}^{-1}_{GFy}}; (23,-6)*{\scriptstyle \Downarrow
F\tilde{\Psi}_{Fy}}; (45,0)*{\cong_{\alpha}}; {\ar^{F \beta_{y}}
(-30,0)*+{Fy}; (0,0)*+{FGFy} }
\endxy
\]
Applying $F$ to the second biadjoint biequivalence axiom and
rewriting, we see that the above pasting diagram is the identity
if and only if the following one is.
\[
\xy {\ar^{F \beta_{y}} (0,0)*+{Fy}; (18,15)*+{FGFy} }; {\ar^{1}
(18,15)*+{FGFy}; (54,15)*+{FGFy} }; {\ar^{\alpha_{Fy}}
(54,15)*+{FGFy}; (72,0)*+{Fy} }; {\ar_{F \beta_{y}} (0,0)*+{Fy};
(18,-15)*+{FGFy} }; {\ar_{1} (18,-15)*+{FGFy}; (54,-15)*+{FGFy} };
{\ar_{\alpha_{Fy}}  (54,-15)*+{FGFy}; (72,0)*+{Fy} }; {\ar_{F
\beta_{GFy}} (18,15)*+{FGFy}; (36,0)*+{FGFGFy} };
{\ar_{\alpha_{FGFy}} (36,0)*+{FGFGFy}; (54,15)*+{FGFy} };
{\ar^{FGF \beta_{y}} (18,-15)*+{FGFy}; (36,0)*+{FGFGFy} };
{\ar^{FG \alpha_{Fy}} (36,0)*+{FGFGFy}; (54,-15)*+{FGFy} };
(18,0)*{\cong_{F \beta}}; (36,9)*{\Downarrow
\tilde{\Phi}^{-1}_{GFy}}; (36,-9)*{\Downarrow FG
\tilde{\Phi}_{y}}; (54,0)*{\cong_{\alpha}}
\endxy
\]

Recall now that $\beta$ was constructed together with an
invertible modification $\delta$ with components $\delta_{y}: F
\beta_{y} \cong \alpha_{Fy}^{\dot}$ such that the composite
\[
\alpha_{Fy} \circ \alpha_{Fy}^{\dot} \stackrel{1 *
\delta_{y}^{-1}}{\Longrightarrow} \alpha_{Fy} \circ F \beta_{y}
\stackrel{\Phi_{y}}{\Longrightarrow} 1_{Fy}
\]
is the counit of the adjoint equivalence $\alpha \dashv_{eq}
\alpha^{\dot}$.  The previous pasting diagram is then the identity
if and only if we pre-compose it with $\delta_{y}^{-1}$ and
post-compose it with $\delta_{y}$.  Using the naturality axiom for
$\alpha^{\dot}$, the modification axiom for $\delta$, and the
equality relating $\delta$ to the counit of the adjoint
equivalence $\alpha \dashv_{eq} \alpha^{\dot}$, the pre- and
post-composed pasting is the identity if and only if the one
displayed below is the identity.
\[
\xy {\ar^{\alpha_{Fy}^{\dot}} (0,0)*+{Fy}; (18,15)*+{FGFy} };
{\ar^{1} (18,15)*+{FGFy}; (54,15)*+{FGFy} }; {\ar^{\alpha_{Fy}}
(54,15)*+{FGFy}; (72,0)*+{Fy} }; {\ar_{\alpha_{Fy}^{\dot}}
(0,0)*+{Fy}; (18,-15)*+{FGFy} }; {\ar_{1} (18,-15)*+{FGFy};
(54,-15)*+{FGFy} }; {\ar_{\alpha_{Fy}}  (54,-15)*+{FGFy};
(72,0)*+{Fy} }; {\ar_{\alpha_{FGFy}^{\dot}} (18,15)*+{FGFy};
(36,0)*+{FGFGFy} }; {\ar_{\alpha_{FGFy}} (36,0)*+{FGFGFy};
(54,15)*+{FGFy} }; {\ar^{FG \alpha_{Fy}^{\dot}} (18,-15)*+{FGFy};
(36,0)*+{FGFGFy} }; {\ar^{FG \alpha_{Fy}} (36,0)*+{FGFGFy};
(54,-15)*+{FGFy} }; (18,0)*{\cong_{\alpha_{F}^{\dot}}};
(36,9)*{\cong_{c^{-1}}}; (36,-9)*{\cong_{c}};
(54,0)*{\cong_{\alpha_{F}}}
\endxy
\]
Here we have written $\cong_{c}$ for the counit isomorphism of the
adjoint equivalence $\alpha \dashv \alpha^{\dot}$, and
$\cong_{c^{-1}}$ for the inverse of the counit.  This diagram is
the identity following from a general lemma on mates that we give
below.
\end{proof}

\lemma
Let $B$ be a bicategory, and $T:B \rightarrow B$ a functor.  Let
$\alpha:T \Rightarrow 1$, $\alpha^{\dot}:1 \Rightarrow T$ be part
of an adjoint equivalence $\alpha \dashv_{eq} \alpha^{\dot}$.
Then for any object $a$ in $B$, the pasting diagram below is the
identity.
\[
\xy {\ar^{\alpha_{a}^{\dot}} (0,0)*+{a}; (18,15)*+{Ta} }; {\ar^{1}
(18,15)*+{Ta}; (54,15)*+{Ta} }; {\ar^{\alpha_{a}} (54,15)*+{Ta};
(72,0)*+{a} }; {\ar_{\alpha_{a}^{\dot}} (0,0)*+{a}; (18,-15)*+{Ta}
}; {\ar_{1} (18,-15)*+{Ta}; (54,-15)*+{Ta} }; {\ar_{\alpha_{a}}
(54,-15)*+{Ta}; (72,0)*+{a} }; {\ar_{\alpha_{Ta}^{\dot}}
(18,15)*+{Ta}; (36,0)*+{T^{2}a} }; {\ar_{\alpha_{Ta}}
(36,0)*+{T^{2}a}; (54,15)*+{Ta} }; {\ar^{T \alpha_{a}^{\dot}}
(18,-15)*+{Ta}; (36,0)*+{T^{2}a} }; {\ar^{T \alpha_{a}}
(36,0)*+{T^{2}a}; (54,-15)*+{Ta} };
(18,0)*{\cong_{\alpha^{\dot}}}; (36,9)*{\cong_{c^{-1}}};
(36,-9)*{\cong_{c}}; (54,0)*{\cong_{\alpha}}
\endxy
\]
\endlemma
\begin{proof}
First, the above pasting diagram is the identity if and only if
the one below is.
\[
\xy {\ar^{\alpha_{a}^{\dot}} (0,0)*+{a}; (18,15)*+{Ta} }; {\ar^{1}
(18,15)*+{Ta}; (54,15)*+{Ta} }; {\ar^{\alpha_{a}} (54,15)*+{Ta};
(72,0)*+{a} }; {\ar^{\alpha_{a}^{\dot}} (0,0)*+{a}; (18,-15)*+{Ta}
}; {\ar_{1} (18,-15)*+{Ta}; (54,-15)*+{Ta} }; {\ar_{\alpha_{a}}
(54,-15)*+{Ta}; (72,0)*+{a} }; {\ar_{\alpha_{Ta}^{\dot}}
(18,15)*+{Ta}; (36,0)*+{T^{2}a} }; {\ar_{\alpha_{Ta}}
(36,0)*+{T^{2}a}; (54,15)*+{Ta} }; {\ar^{T \alpha_{a}^{\dot}}
(18,-15)*+{Ta}; (36,0)*+{T^{2}a} }; {\ar^{T \alpha_{a}}
(36,0)*+{T^{2}a}; (54,-15)*+{Ta} };
(18,0)*{\cong_{\alpha^{\dot}}}; (36,9)*{\cong_{c^{-1}}};
(36,-9)*{\cong_{c}}; (54,0)*{\cong_{\alpha}};
{\ar^{\alpha_{a}}
(-18,0)*+{Ta}; (0,0)*+{a} }; {\ar@/^1.5pc/^{1} (-18,0)*+{Ta};
(18,15)*+{Ta} }; {\ar@/_1.5pc/_{1} (-18,0)*+{Ta}; (18,-15)*+{Ta}
}; (0,6)*{\cong_{u}}; (0,-6)*{\cong_{u^{-1}}}
\endxy
\]
Here we have written $\cong_{u}$ to denote a unit isomorphism, and
$\cong_{u^{-1}}$ the inverse of a unit isomorphism.

It is now a simple calculation using the adjoint equivalence
axioms for transformations to show that the pasting below (modulo unit
isomorphisms to alter the 1-cell source and target, which we
ignore for now but record the presence of for later)
\[
\xy
{\ar^{\alpha_{a}} (0,0)*+{Ta}; (25,0)*+{a} };
{\ar^{\alpha_{a}^{\dot}} (25,0)*+{a}; (50,0)*+{Ta} };
{\ar^{1} (50,0)*+{Ta}; (75,-14)*+{Ta} };
{\ar_{1} (0,0)*+{Ta}; (25,-14)*+{Ta} };
{\ar_{T \alpha_{a}^{\dot}} (25,-14)*+{Ta}; (50,-14)*+{T^{2} a} };
{\ar_{\alpha_{Ta}} (50,-14)*+{T^{2} a}; (75,-14)*+{Ta} };
{\ar^{\alpha_{a}^{\dot}} (25,0)*+{a}; (25,-14)*+{Ta} };
{\ar_{\alpha_{Ta}^{\dot}} (50,0)*+{Ta}; (50,-14)*+{T^{2} a} };
(37.5,-7)*{\cong_{\alpha^{\dot}}}; (17,-5)*{\cong_{u^{-1}}}; (59,-9)*{\cong_{c^{-1}}}
\endxy
\]
is equal to
the naturality square below.
\[ \xy
\endxy
\]
Thus we have shown that the pasting diagram in the previous paragraph is equal (modulo units) to the pasting diagram displayed below.
\[
\xy {\ar^{\alpha_{a}} (0,0)*+{Ta}; (15,9)*+{a} };
{\ar^{\alpha_{a}^{\dot}} (15,9)*+{a}; (30,15)*+{Ta} };
{\ar@/^3pc/^{1} (0,0)*+{Ta}; (30,15)*+{Ta} };
{\ar_{T\alpha_{a}^{\dot}} (0,0)*+{Ta}; (30,0)*+{T^{2}a} };
{\ar_{\alpha_{Ta}} (30,0)*+{T^{2}a}; (30,15)*+{Ta} };
{\ar^{\alpha_{a}}  (30,15)*+{Ta}; (60,15)*+{a} }; {\ar_{T
\alpha_{a}} (30,0)*+{T^{2}a}; (60,0)*+{Ta} }; {\ar_{\alpha_{a}}
(60,0)*+{Ta}; (60,15)*+{a} }; {\ar@/_2pc/_{1} (0,0)*+{Ta};
(60,0)*+{Ta} }; (12,15)*{\cong_{u}}; (18,6)*{\cong_{\alpha}};
(45,7.5)*{\cong_{\alpha}}; (30,-5)*{\cong_{c}}
\endxy
\]
Using the triangle identities and the naturality axioms for
$\alpha$, this is now equal to the pasting below (once again
modulo units).
\[
\xy {\ar^{\alpha_{a}} (0,0)*+{Ta}; (0,15)*+{a} }; {\ar^{\alpha_{a}
\alpha_{a}^{\dot}} (0,15)*+{a}; (30,15)*+{a} };
{\ar^{T(\alpha_{a}\alpha_{a}^{\dot})} (0,0)*+{Ta}; (30,0)*+{Ta} };
{\ar_{\alpha_{a}} (30,0)*+{Ta}; (30,15)*+{a} }; {\ar@/^3pc/^{1}
(0,15)*+{a}; (30,15)*+{a} }; {\ar@/_2pc/_{1} (0,0)*+{Ta};
(30,0)*+{Ta} }; (15,7.5)*{\cong_{\alpha}};
(15,22)*{\cong_{c^{-1}}}; (15,-4)*{\cong_{c}}
\endxy
\]
By the naturality axioms for $\alpha$ and coherence for functors,
this is equal to a composite of left and right unit isomorphisms
so the original diagram is a composite of coherence isomorphisms,
hence is the identity by coherence for bicategories.
\end{proof}

\section{Biequivalences in general tricategories}

This section will establish the general result that every
biequivalence in a tricategory $T$ is part of a biadjoint
biequivalence in $T$.  Our proof will proceed largely as did the
general case for equivalences in bicategories by first examining
the case of functor tricategories and then using a Yoneda
embedding.  Since a Yoneda embedding is only known for cubical
tricategories rather than the general case, the proof for
tricategories is slightly longer although essentially the same.
We refer the reader to \cite{gps} or \cite{gur} for the relevant tricategorial
results.


Recall that if $S$ is any tricategory and $T$ is a
$\textbf{Gray}$-category, then there is a $\textbf{Gray}$-category
$\textbf{Tricat}(S,T)$ with objects functors $S \rightarrow T$,
1-cells transformations, 2-cells modifications, and 3-cells
perturbations.


\proposition
Let $S$ be any tricategory and $T$ be a $\textbf{Gray}$-category.
Assume that every biequivalence 1-cell $f$ in $T$ is part of a
biadjoint biequivalence $f \dashv_{bieq} g$.  Then every
biequivalence 1-cell $\alpha$ in $\textbf{Tricat}(S,T)$ is part of
a biadjoint biequivalence $\alpha \dashv_{bieq} \beta$.
\endproposition

\remark
We have written the proof of this proposition out so that it should be obvious to the reader that it remains true when $T$ is merely a tricategory and not a $\textbf{Gray}$-category.  By this we mean the following:  the assumption that $T$ is a $\textbf{Gray}$-category is only present to use the results of \cite{gur} to give a concrete construction of the tricategory $\textbf{Tricat}(S,T)$.  Using coherence for tricategories, it is possible to construct a tricategory $\textbf{Tricat}(S,T)$ when $T$ is any tricategory, and then the proof below applies verbatim to the analogous proposition.
\endremark

\begin{proof}
Let $\alpha:F \Rightarrow G$ be a biequivalence in
$\textbf{Tricat}(S,T)$.  Then $\alpha_{a}:Fa \rightarrow Ga$ is a
biequivalence in $T$ for every object $a$, so we choose biadjoint
biequivalences $\alpha_{a} \dashv_{bieq} \beta_{a}$ for every
object $a$ of $S$.  To complete the proof, we must do the
following:
\begin{enumerate}
\item equip the components $\beta_{a}$ with the structure of a
transformation;
\item equip the componentwise adjoint equivalences
\[
\begin{array}{c}
\bs{\varepsilon}_{a}:\alpha_{a} \beta_{a} \dashv_{eq} 1_{Ga},
\quad \bs{\eta}_{b}: 1_{Fb} \dashv_{eq} \beta_{b} \alpha_{b}
\end{array}
\]
with the structure of adjoint equivalences in the
hom-bicategories $\textbf{Tricat}(S,T)(G,G)$,
$\textbf{Tricat}(S,T)(F,F)$; and
\item check that
$\Phi, \Psi$ are perturbations.
\end{enumerate}
Since equations between perturbations are checked componentwise,
the fact we have biadjoint biequivalences $\alpha_{a}
\dashv_{bieq} \beta_{a}$ will then imply that there is a global
biadjoint biequivalence $\alpha \dashv_{bieq} \beta$ in
$\textbf{Tricat}(S,T)$.

We begin by defining a transformation $\beta:G \Rightarrow F$ with
components given by these $\beta_{a}$.  Since we have already
given the components on objects, there are three pieces of data
left to define.  The first is an adjoint equivalence
\[
\xy {\ar^{G} (0,0)*+{S(a,b)}; (30,0)*+{T(Ga, Gb)} }; {\ar^{T(1,
\beta_{b})} (30,0)*+{T(Ga, Gb)}; (30,-15)*+{T(Ga, Fb)} }; {\ar_{F}
(0,0)*+{S(a,b)}; (0,-15)*+{T(Fa, Fb)} }; {\ar_{T(\beta_{a},1)}
(0,-15)*+{T(Fa, Fb)}; (30,-15)*+{T(Ga, Fb)} };
{\ar@{=>}_{\bs{\beta}} (22,-4)*{}; (8,-11)*{} }
\endxy
\]
in the bicategory $\textbf{Bicat}\Big( S(a,b), T(Ga,Fb) \Big)$.
We write down the component $\beta_{f}$ of $\beta$ at an object
$f$ and leave it to the reader to construct the rest of the
adjoint equivalence in the obvious fashion.
\[
\begin{array}{rcl}
\beta_{b} \otimes Gf & \stackrel{1 \otimes r^{\dot}}{\longrightarrow} & \beta_{b} \otimes (Gf \otimes 1) \\
& \stackrel{1 \otimes (1 \otimes \varepsilon_{a}^{\dot})}{\longrightarrow} & \beta_{b} \otimes (Gf \otimes (\alpha_{a} \otimes \beta_{a})) \\
& \stackrel{1 \otimes a^{\dot}}{\longrightarrow} & \beta_{b} \otimes ((Gf \otimes \alpha_{a}) \otimes \beta_{a}) \\
& \stackrel{1 \otimes (\alpha_{f}^{\dot} \otimes 1)}{\longrightarrow} & \beta_{b} \otimes ((\alpha_{b} \otimes Ff) \otimes \beta_{a}) \\
& \stackrel{a^{\dot}}{\longrightarrow} & (\beta_{b} \otimes (\alpha_{b} \otimes Ff)) \otimes \beta_{a} \\
& \stackrel{a^{\dot} \otimes 1}{\longrightarrow} & ((\beta_{b} \otimes \alpha_{b})\otimes Ff) \otimes \beta_{a} \\
& \stackrel{\eta_{b}^{\dot} \otimes 1}{\longrightarrow} & (1 \otimes Ff) \otimes \beta_{a} \\
& \stackrel{l \otimes 1}{\longrightarrow} & Ff \otimes \beta_{a}
\end{array}
\]
We have written this out as if if were a 2-cell an arbitrary tricategory, not necessarily a $\mb{Gray}$-category.  In the case that $T$ is $\mb{Gray}$, this cell is as below.
\[
\beta_{b} Gf  \stackrel{11 \epsilon^{\dot}_{a}}{\longrightarrow}
\beta_{a} Gf \alpha_{a} \beta_{a}  \stackrel{1 \alpha_{f}^{\dot} 1}{\longrightarrow} \beta_{b} \alpha_{b} Ff \beta_{a} \stackrel{\eta_{b}^{\dot} 1}{\longrightarrow} Ff \beta_{a}
\]

We now must produce a pair of invertible modifications $\Pi, M$ to
complete the definition of the data for the transformation
$\beta$.
The component of the modification $M^{\beta}$ at the object $a$ is
given by the isomorphism shown below where the unmarked
isomorphisms are unique by coherence and the two marked 3-cells
are both appropriate mates.
\[
\xy {\ar^{\scriptstyle 1 \iota^{G}} (0,0)*+{\scriptstyle
\beta_{a}}; (25,0)*+{\scriptstyle \beta_{a} GI} };
{\ar^{\scriptstyle 11 \varepsilon^{\dot}}  (25,0)*+{\scriptstyle
\beta_{a} GI}; (50,0)*+{\scriptstyle \beta_{a} GI \alpha_{a}
\beta_{a}} };
{\ar^{\scriptstyle 1 \alpha_{I}^{\dot} 1}
(50,0)*+{\scriptstyle \beta_{a} GI \alpha_{a} \beta_{a}};
(75,0)*+{\scriptstyle \beta_{a} \alpha_{a} FI \beta_{a}} };
{\ar^{\scriptstyle \eta^{\dot} 11} (75,0)*+{\scriptstyle \beta_{a}
\alpha_{a} FI \beta_{a}}; (75,-28)*+{\scriptstyle FI \beta_{a}} };
{\ar_{\scriptstyle 1 \varepsilon^{\dot}}  (0,0)*+{\scriptstyle
\beta_{a}}; (25,-24)*+{\scriptstyle \beta_{a} \alpha_{a}
\beta_{a}} };
{\ar_{\scriptstyle 1 \iota^{G} 11}
(25,-24)*+{\scriptstyle \beta_{a} \alpha_{a} \beta_{a}};
(50,0)*+{\scriptstyle \beta_{a} GI \alpha_{a} \beta_{a}} };
{\ar@/_2pc/_{\scriptstyle 11 \iota^{F} 1} (25,-24)*+{\scriptstyle
\beta_{a} \alpha_{a} \beta_{a}}; (75,0)*+{\scriptstyle \beta_{a}
\alpha_{a} FI \beta_{a}} }; {\ar^{\scriptstyle \eta^{\dot} 1}
(25,-24)*+{\scriptstyle \beta_{a} \alpha_{a} \beta_{a}};
(25,-32)*+{\scriptstyle \beta_{a}} }; {\ar@/_2pc/_{\scriptstyle 1}
(0,0)*+{\scriptstyle \beta_{a}}; (25,-32)*+{\scriptstyle
\beta_{a}} }; {\ar@/_2pc/_{\scriptstyle \iota^{F} 1}
(25,-32)*+{\scriptstyle \beta_{a}}; (75,-28)*+{\scriptstyle FI
\beta_{a}} }; (12,-20)*{\scriptstyle \Downarrow \Psi^{-1}};
(25,-7)*{\scriptstyle \cong}; (54,-7)*{\scriptstyle \Downarrow
1(M_{a}^{\alpha})^{-1}1}; (50,-29)*{\scriptstyle \cong}
\endxy
\]
The component of the modification $\Pi^{\beta}$ at the composable
pair $(g,f)$ is given by the pasting below, once again following
the same conventions.  To conserve space, we omit the subscripts for the components of $\alpha$ and $\beta$ given that they can be deduced from the other 1-cells in any given term.
\[
\xy {\ar^{\scriptstyle 11\epsilon^{\dot} 1} (0,0)*+{\scriptstyle
\beta GgGf}; (25,0)*+{\scriptstyle \beta Gg \alpha
\beta Gf} }; {\ar^{\scriptstyle 1 \alpha^{\dot}}
(25,0)*+{\scriptstyle \beta Gg \alpha \beta Gf};
(60,0)*+{\scriptstyle \beta \alpha Fg \beta Gf} };
{\ar^{\scriptstyle \eta^{\dot}111} (60,0)*+{\scriptstyle
\beta \alpha Fg \beta Gf}; (85,0)*+{\scriptstyle Fg
\beta Gf} }; {\ar^{\scriptstyle 111\epsilon^{\dot}}
(85,0)*+{\scriptstyle Fg \beta Gf}; (85,-15)*+{\scriptstyle
Fg \beta Gf \alpha \beta} }; {\ar^{\scriptstyle 11
\alpha^{\dot}1} (85,-15)*+{\scriptstyle Fg \beta Gf
\alpha \beta}; (85,-30)*+{\scriptstyle Fg \beta
\alpha Ff \beta} }; {\ar^{\scriptstyle 1 \eta^{\dot}11}
(85,-30)*+{\scriptstyle Fg \beta \alpha Ff \beta};
(85,-45)*+{\scriptstyle Fg Ff \beta} }; {\ar_{\scriptstyle
111 \epsilon^{\dot}} (0,0)*+{\scriptstyle \beta GgGf};
(0,-15)*+{\scriptstyle \beta Gg Gf \alpha \beta} };
{\ar^{\scriptstyle 11 \epsilon^{\dot}111}  (0,-15)*+{\scriptstyle
\beta Gg Gf \alpha \beta}; (25,-15)*+{\scriptstyle
\beta Gg \alpha \beta Gf \alpha \beta} };
{\ar^{\scriptstyle 1 \alpha^{\dot}1111}
(25,-15)*+{\scriptstyle \beta Gg \alpha \beta Gf
\alpha \beta}; (60,-15)*+{\scriptstyle \beta
\alpha Fg \beta Gf \alpha \beta} };
{\ar^{\scriptstyle \eta^{\dot}11111} (60,-15)*+{\scriptstyle
\beta \alpha Fg \beta Gf \alpha \beta};
(85,-15)*+{\scriptstyle
Fg \beta Gf \alpha \beta} };
{\ar^{\scriptstyle 11 \alpha^{\dot}1} (0,-15)*+{\scriptstyle
\beta Gg Gf \alpha \beta}; (0,-30)*+{\scriptstyle
\beta Gg \alpha Ff \beta} }; {\ar^{\scriptstyle
11\epsilon^{\dot}111} (0,-30)*+{\scriptstyle \beta Gg
\alpha Ff \beta}; (25,-30)*+{\scriptstyle \beta Gg
\alpha \beta \alpha Ff \beta} };
{\ar^{\scriptstyle 1 \alpha^{\dot}1111}
(25,-30)*+{\scriptstyle \beta Gg \alpha \beta
\alpha Ff \beta}; (60,-30)*+{\scriptstyle \beta
\alpha Fg \beta \alpha Ff \beta} };
{\ar^{\scriptstyle \eta^{\dot}11111} (60,-30)*+{\scriptstyle
\beta \alpha Fg \beta \alpha Ff \beta};
(85,-30)*+{\scriptstyle Fg \beta \alpha Ff \beta} };
{\ar^{\scriptstyle 111\eta^{\dot}11} (25,-30)*+{\scriptstyle
\beta Gg \alpha \beta \alpha Ff \beta};
(25,-45)*+{\scriptstyle \beta Gg \alpha Ff \beta} };
{\ar_{\scriptstyle 1 \alpha^{\dot}11} (25,-45)*+{\scriptstyle
\beta Gg \alpha Ff \beta}; (60,-45)*+{\scriptstyle
\beta \alpha Fg Ff \beta} }; {\ar_{\scriptstyle
\eta^{\dot}111} (60,-45)*+{\scriptstyle \beta \alpha Fg
Ff \beta}; (85,-45)*+{\scriptstyle Fg Ff \beta} };
{\ar@/_2pc/_{\scriptstyle 1} (0,-30)*+{\scriptstyle \beta Gg
\alpha Ff \beta}; (25,-45)*+{\scriptstyle \beta Gg
\alpha Ff \beta} }; {\ar_{\scriptstyle 1\chi^{G}}
(0,0)*+{\scriptstyle \beta GgGf}; (-20,-10)*+{\scriptstyle
\beta G(gf)} }; {\ar_{\scriptstyle 11 \epsilon^{\dot}}
(-20,-10)*+{\scriptstyle \beta G(gf)};
(-20,-25)*+{\scriptstyle \beta G(gf) \alpha \beta} };
{\ar^{\scriptstyle 1 \chi^{G}11} (0,-15)*+{\scriptstyle \beta
Gg Gf \alpha \beta}; (-20,-25)*+{\scriptstyle \beta
G(gf) \alpha \beta} }; {\ar@/_4pc/_{\scriptstyle 1
\alpha^{\dot} 1} (-20,-25)*+{\scriptstyle \beta G(gf)
\alpha \beta}; (60,-60)*+{\scriptstyle \beta
\alpha F(gf) \beta} }; {\ar_{\scriptstyle \eta^{\dot}11}
(60,-60)*+{\scriptstyle \beta \alpha F(gf) \beta};
(85,-60)*+{\scriptstyle F(gf) \beta} }; {\ar_{\scriptstyle 11
\chi^{F}1} (60,-45)*+{\scriptstyle \beta \alpha Fg Ff
\beta}; (60,-60)*+{\scriptstyle \beta \alpha F(gf)
\beta} }; {\ar^{\scriptstyle \chi^{F} 1}
(85,-45)*+{\scriptstyle Fg Ff \beta}; (85,-60)*+{\scriptstyle
F(gf) \beta} }; (37.5,-7.5)*{\scriptstyle \cong};
(37.5,-22.5)*{\scriptstyle \cong}; (50,-37.5)*{\scriptstyle
\cong}; (-10,-12)*{\scriptstyle \cong};
(72.5,-52.5)*{\scriptstyle \cong}; (11,-36)*{\scriptstyle
\Downarrow 11\Psi^{-1}11}; (15,-52)*{\scriptstyle \Downarrow 1
(\Pi_{gf}^{\alpha})^{-1}1}
\endxy
\]
All of the unmarked isomorphisms are isomorphisms of the form
\[
(\beta \otimes 1) \circ (1 \otimes \alpha) \cong (1 \otimes \alpha) \circ (\beta \otimes 1)
\]
arising from the functoriality of the horizontal composition
\[
\otimes:T(y,z) \times T(x,y) \rightarrow T(x,z).
\]
The transformation axioms for $\beta$ then follow from those
for $\alpha$, thus completing the construction of a weak inverse
for $\alpha$.

The next step is to construct adjoint equivalences
\[
\begin{array}{c}
\bs{\varepsilon}:\alpha \beta \dashv_{eq} 1_{G}, \quad \bs{\eta}:
1_{F} \dashv_{eq} \beta \alpha.
\end{array}
\]
We already have the adjoint equivalences
\[
\begin{array}{c}
\bs{\varepsilon}_{a}:\alpha_{a} \beta_{a} \dashv_{eq} 1_{Ga},
\quad \bs{\eta}_{b}: 1_{Fb} \dashv_{eq} \beta_{b} \alpha_{b}
\end{array}
\]
on components in the tricategory $T$, we need only lift these to
adjoint equivalences in the hom-bicategories
$\textbf{Tricat}(S,T)(G,G), \textbf{Tricat}(S,T)(F,F)$,
respectively.  Thus we need to equip the collection
$\varepsilon_{a}: \alpha_{a} \beta_{a} \Rightarrow I_{Ga}$ with
the structure of a trimodification.  The 2-cell $\varepsilon_{a}$
is the 2-cell in $T$ of the same name, so we must give the
invertible 3-cell displayed below.
\[
\xy {\ar^{(\alpha\beta)_{f}} (0,0)*+{\alpha_{b}\beta_{b}Gf};
(20,12)*+{Gf\alpha_{a}\beta_{a}} }; {\ar^{1\varepsilon_{a}}
(20,12)*+{Gf\alpha_{a}\beta_{a}}; (40,0)*+{Gf 1_{Ga}} };
{\ar_{\varepsilon_{b}1} (0,0)*+{\alpha_{b}\beta_{b}Gf};
(20,-12)*+{1_{Gb}Gf} }; {\ar_{(1_{G})_{f}}  (20,-12)*+{1_{Gb}Gf};
(40,0)*+{Gf 1_{Ga}} }; {\ar@{=>}^{\varepsilon_{f}} (20,7)*{};
(20,-7)*{} }
\endxy
\]
The component $(1_{G})_{f}$ is given by $r^{\dot} \otimes l$, so
once again using coherence we write this as the identity and can
define $\varepsilon_{f}$ as shown below.

\[
\xy
{\ar^{\scriptstyle 111 \epsilon^{\dot}} (0,0)*+{\scriptstyle \alpha_{b} \beta_{b} Gf};
(0,15)*+{\scriptstyle \alpha_{b} \beta_{b} Gf \alpha_{a} \beta_{a}} };
{\ar^{\scriptstyle 11 \alpha_{f}^{\dot}1} (0,15)*+{\scriptstyle \alpha_{b} \beta_{b} Gf
\alpha_{a} \beta_{a}}; (15,30)*+{\scriptstyle \alpha_{b} \beta_{b} \alpha_{b}
Ff \beta_{a}} };
{\ar@/^1pc/^{\scriptstyle 1 \eta^{\dot} 11} (15,30)*+{\scriptstyle \alpha_{b} \beta_{b}
\alpha_{b} Ff \beta_{a}}; (65,30)*+{\scriptstyle \alpha_{b} Ff \beta_{a}} };
{\ar^{\scriptstyle \alpha_{f} 1} (65,30)*+{\scriptstyle \alpha_{b} Ff \beta_{a}};
(80,15)*+{\scriptstyle Gf \alpha_{a} \beta_{a}} };
{\ar^{\scriptstyle 1 \epsilon} (80,15)*+{\scriptstyle Gf \alpha_{a} \beta_{a}};
(80,0)*+{\scriptstyle Gf} };
{\ar_{\scriptstyle 1} (0,0)*+{\scriptstyle \alpha_{b} \beta_{b} Gf}; (40,0)*+{\scriptstyle \alpha_{b} \beta_{b} Gf} };
{\ar_{\scriptstyle \epsilon 1} (40,0)*+{\scriptstyle \alpha_{b} \beta_{b} Gf}; (80,0)*+{\scriptstyle Gf}
};
{\ar_{\scriptstyle 1} (0,15)*+{\scriptstyle \alpha_{b} \beta_{b} Gf \alpha_{a} \beta_{a}}; (40,12)*+{\scriptstyle \alpha_{b} \beta_{b} Gf \alpha_{a} \beta_{a}} };
{\ar_{\scriptstyle 111 \epsilon} (40,12)*+{\scriptstyle \alpha_{b} \beta_{b} Gf \alpha_{a} \beta_{a}}; (40,0)*+{\scriptstyle \alpha_{b} \beta_{b} Gf} };
{\ar^{\scriptstyle \epsilon^{\dot}11} (80,15)*+{\scriptstyle Gf \alpha_{a} \beta_{a}};
(40,12)*+{\scriptstyle \alpha_{b} \beta_{b} Gf \alpha_{a} \beta_{a}} };
{\ar_{\scriptstyle 1} (15,30)*+{\scriptstyle \alpha_{b} \beta_{b} \alpha_{b} Ff \beta_{a}}; (40,27)*+{\scriptstyle \alpha_{b} \beta_{b} \alpha_{b} Ff \beta_{a}} };
{\ar^{\scriptstyle \epsilon^{\dot}111} (65,30)*+{\scriptstyle \alpha_{b} Ff \beta_{a}};
(40,27)*+{\scriptstyle \alpha_{b} \beta_{b} \alpha_{b} Ff \beta_{a}} };
{\ar_{\scriptstyle 11 \alpha_{f} 1} (40,27)*+{\scriptstyle \alpha_{b} \beta_{b} \alpha_{b}Ff \beta_{a}}; (40,12)*+{\scriptstyle \alpha_{b} \beta_{b} Gf \alpha_{a}
\beta_{a}} };
(20,7)*{\scriptstyle \cong}; (60,7)*{\scriptstyle \cong}; (20,20)*{\scriptstyle \cong}; (60,20)*{\scriptstyle \cong};
(40,32)*{\scriptstyle \Leftarrow \Psi^{-1}11}
\endxy
\]
The unmarked isomorphisms in this pasting diagram are (up to coherence) counits for adjoint equivalences for the two
left squares, functoriality of the tensor for the top right square, and a counit together with functoriality of the tensor for the bottom right square.  There are now two trimodification axioms to check, one for composition and one for units.  Both of these axioms follows from coherence, the definitions of $\Pi^{\beta}$ and $M^{\beta}$, and the definitions of $\Pi^{\beta \alpha}$ and $M^{\beta \alpha}$; they do not require the biadjunction axioms for $\alpha_{a} \dashv_{bieq} \beta_{a}$.

To complete the proof that we have an adjoint equivalence
\[
\mb{\varepsilon}: \alpha \beta \dashv_{eq} 1_{G}
\]
in $\mb{Tricat}(S,T)(G,G)$, we need only note that a modification $m: \theta \rightarrow \phi$ in this bicategory is an equivalence if and only if each component $m_{x}:\theta_{x} \Rightarrow \phi_{x}$ is an equivalence in the appropriate hom-bicategory of $T$.  In this case, the modification $m$ is $\varepsilon$, and each component $\varepsilon_{a}$ is an equivalence in $T(Ga,Ga)$ by construction.  We note in passing that the components of the specified pseudo-inverse for $\varepsilon$ can be taken to be the 2-cells $\varepsilon_{a}^{\dot}$ that appear in the pointwise adjoint equivalences
\[
\mb{\varepsilon}: \alpha_{a} \beta_{a} \dashv_{eq} 1_{Ga}.
\]

The construction of the adjoint equivalence
\[
\mb{\eta}: 1_{F} \dashv_{eq} \beta \alpha
\]
follows exactly the same pattern as for $\mb{\varepsilon}$, so we omit most of the details here.  In order to check that $\Phi$ and $\Psi$ are perturbations, we will need to know the components
\[
\eta_{f}: 1 \eta_{a} \circ (1_{F})_{f} \Rrightarrow (\beta \alpha)_{f} \circ \eta_{b}1
\]
for each 1-cell $f:a \rightarrow b$.  Using coherence, we write $\eta_{f}$ as the pasting below.
\[
\xy
{\ar^{\scriptstyle 1 \eta} (0,0)*+{\scriptstyle Ff}; (40,0)*+{\scriptstyle Ff \beta \alpha} };
{\ar^{\scriptstyle 1} (40,0)*+{\scriptstyle Ff \beta \alpha}; (80,0)*+{\scriptstyle Ff \beta \alpha} };
{\ar_{\scriptstyle \eta 1} (0,0)*+{\scriptstyle Ff}; (0,-15)*+{\scriptstyle \beta \alpha Ff} };
{\ar_{\scriptstyle 1 \alpha_{f}} (0,-15)*+{\scriptstyle \beta \alpha Ff}; (15,-30)*+{\scriptstyle \beta Gf \alpha} };
{\ar@/_1pc/_{\scriptstyle 11\varepsilon^{\dot}1} (15,-30)*+{\scriptstyle \beta Gf \alpha}; (65,-30)*+{\scriptstyle \beta Gf \alpha \beta \alpha} };
{\ar_{\scriptstyle 1 \alpha_{f}^{\dot} 11} (65,-30)*+{\scriptstyle \beta Gf \alpha \beta \alpha}; (80,-15)*+{\scriptstyle \beta \alpha Ff \beta \alpha} };
{\ar_{\scriptstyle \eta^{\dot}111} (80,-15)*+{\scriptstyle \beta \alpha Ff \beta \alpha}; (80,0)*+{\scriptstyle Ff \beta \alpha} };
{\ar^{\scriptstyle \eta 111} (40,0)*+{\scriptstyle Ff \beta \alpha}; (40,-12)*+{\scriptstyle \beta \alpha Ff \beta \alpha} };
{\ar^{\scriptstyle 1} (40,-12)*+{\scriptstyle \beta \alpha Ff \beta \alpha}; (80,-15)*+{\scriptstyle \beta \alpha Ff \beta \alpha} };
{\ar_{\scriptstyle 111 \eta^{\dot}} (40,-12)*+{\scriptstyle \beta \alpha Ff \beta \alpha}; (0,-15)*+{\scriptstyle \beta \alpha Ff} };
{\ar^{\scriptstyle 1 \alpha_{f} 11} (40,-12)*+{\scriptstyle \beta \alpha Ff \beta \alpha}; (40,-27)*+{\scriptstyle \beta Gf \alpha \beta \alpha} };
{\ar^{\scriptstyle 1} (40,-27)*+{\scriptstyle \beta Gf \alpha \beta \alpha}; (65,-30)*+{\scriptstyle \beta Gf \alpha \beta \alpha} };
{\ar_{\scriptstyle 111 \eta^{\dot}} (40,-27)*+{\scriptstyle \beta Gf \alpha \beta \alpha}; (15,-30)*+{\scriptstyle \beta Gf \alpha} };
(20,-7)*{\scriptstyle \cong}; (20,-20)*{\scriptstyle \cong}; (60,-7)*{\scriptstyle \cong}; (60,-20)*{\scriptstyle \cong};
(40,-32)*{\scriptstyle \Leftarrow 11 \Psi}
\endxy
\]
The isomorphisms are all obtained from units of adjunctions and functoriality of the tensor.

Finally, we must check that the cells $\Phi_{a}, \Psi_{a}$
constitute a pair of invertible perturbations.  This involves checking a single axiom for each 1-cell $f:a \rightarrow b$ which we leave to the reader, but it follows from coherence and the biadjunction axioms for both $\alpha_{a} \dashv_{bieq} \beta_{a}$ and $\alpha_{b} \dashv_{bieq} \beta_{b}$.
\end{proof}


To prove our main result, we require two lemmas.

\lemma
Let $f:a \rightarrow b$ and $g,h:b \rightarrow a$ be 1-cells in a tricategory $T$.  If both $gf$ is equivalent to $1_{a}$ in $T(a,a)$ and $fh$ is equivalent to $1_{b}$ in $T(b,b)$, then $h$ is equivalent to $g$ in $T(b,a)$.  In particular, if $f$ is a biequivalence in $T$ and both $g_{1},g_{2}$ are 1-cells such that $fg_{i} \simeq 1$ and $g_{i}f \simeq 1$ for $i=1,2$, then $g_{1} \simeq g_{2}$.
\endlemma
\begin{proof}
This follows in the standard way by considering the 1-cell $g_{1}fg_{2}$ in $T(b,a)$.
\end{proof}

\lemma
Let $S,T$ be tricategories, and assume that every biequivalence
$f$ in $T$ is part of a biadjoint biequivalence $f \dashv_{bieq}
g$.  Assume that $F:S \rightarrow T$ is 2-locally an equivalence,
i.e., every functor
\[
F:S(a,b)(f,g) \rightarrow T(Fa,Fb)(Ff,Fg)
\]
is an equivalence of categories.  Then every biequivalence $h$ in
$S$ is part of a biadjoint biequivalence.
\endlemma
\begin{proof}
This lemma follows in the same way that Lemma \ref{functorbicats} did; the proof requires using the previous lemma in exactly the same way that Lemma \ref{functorbicats} required the uniqueness of weak inverses in a bicategory.
\end{proof}

We now present our main result.

\theorem\label{biadjbieqintricats}
Let $T$ be a tricategory, and let $f$ be a biequivalence
in $T$.  Then there is a 1-cell $g$ in $T$ and a biadjoint
biequivalence $f \dashv_{bieq} g$.
\endtheorem
\begin{proof}
First, recall that for every tricategory $T$ there is a cubical
tricategory (i.e., a tricategory in which the hom-bicategories are
2-categories and the unit and composition functors are cubical)
$\textrm{st}(T)$ and a triequivalence $T \rightarrow \textrm{st}(T)$.  Since triequivalences satisfy the hypotheses of the lemma
above, we are left proving the theorem in the case that $T$ is a
cubical tricategory.  If $T$ is cubical, then it has a Yoneda
embedding $T \hookrightarrow \textbf{Tricat}(T^{\textrm{op}},
\textbf{Gray})$ which also satisfies the hypotheses of the above
lemma.  Thus if we show that every biequivalence in
$\textbf{Gray}$ is part of a biadjoint biequivalence, we will have
proven the theorem for arbitrary $T$.  But the inclusion
$\textbf{Gray} \hookrightarrow \textbf{Bicat}$ satisfies the
hypotheses of the lemma and we have already proven the claim
directly for $\textbf{Bicat}$.
\end{proof}

\section{Application: lifting monoidal structures}
This section will show how to lift monoidal structures on
bicategories using the results of Section 4.  We refer the reader to \cite{ds} for the definitions of the morphisms and higher cells between monoidal bicategories in the case of $\mb{Gray}$-monoids, and \cite{gps} or \cite{gur} for the definitions of higher cells between tricategories from which these are derived.  From this point forward, most calculations will only be described as they are generally straightforward but the pastings used can be very large.

\theorem
Let $B$ be a monoidal bicategory and let $C$ be any bicategory.
If $F:C \rightarrow B$ is a biequivalence from $C$ to the
underlying bicategory of $B$, then $C$ can be equipped with the
structure of a monoidal bicategory and $F$ can be compatibly
equipped with the structure of a monoidal functor.
\endtheorem
\begin{proof}
First choose a biadjoint biequivalence $F \dashv_{bieq} G$ in
$\textbf{Bicat}$.  We now define a monoidal structure on $C$ as
follows.  The tensor $\boxtimes: C \times C \rightarrow C$ is the
composite
\[
C \times C \stackrel{F \times F}{\longrightarrow} B \times B
\stackrel{\otimes}{\rightarrow} B \stackrel{G}{\rightarrow} C.
\]
The unit $I^{C}:* \rightarrow C$ is the composite
\[
* \stackrel{I^{B}}{\rightarrow} B \stackrel{G}{\rightarrow} C.
\]

The associativity adjoint equivalence is given by the pasting
below.  Here we have marked identity adjoint equivalences with
equal signs.
\[
\xy {\ar^{(F \times F) \times 1} (0,0)*+{(C \times C) \times C};
(33.3,0)*+{(B \times B) \times C} }; {\ar^{\otimes \times 1}
(33.3,0)*+{(B \times B) \times C}; (66.7,0)*+{B \times C} };
{\ar^{G \times 1} (66.7,0)*+{B \times C}; (100,0)*+{C \times C} };
{\ar^{F \times F} (100,0)*+{C \times C}; (100,-15)*+{B \times B}
}; {\ar^{\otimes} (100,-15)*+{B \times B}; (100,-30)*+{B} };
{\ar^{G} (100,-30)*+{B}; (100,-45)*+{C} }; {\ar_{a} (0,0)*+{(C
\times C) \times C}; (0,-15)*+{C \times (C \times C)} }; {\ar_{1
\times (F \times F)} (0,-15)*+{C \times (C \times C)}; (0,-30)*+{C
\times (B \times B)} }; {\ar_{1 \times \otimes} (0,-30)*+{C \times
(B \times B)}; (0,-45)*+{C \times B} }; {\ar_{1 \times G}
(0,-45)*+{C \times B}; (25,-50)*+{C \times C} }; {\ar_{F \times F}
(25,-50)*+{C \times C}; (50,-50)*+{B \times B} }; {\ar_{\otimes}
(50,-50)*+{B \times B}; (75,-50)*+{B} }; {\ar_{G} (75,-50)*+{B};
(100,-45)*+{C} }; {\ar^{1 \times F} (33.3,0)*+{(B \times B) \times
C}; (50,-15)*+{(B \times B) \times B} }; {\ar_{1 \times F}
(66.7,0)*+{B \times C}; (100,-15)*+{B \times B} }; {\ar_{(F \times
F) \times F} (0,0)*+{(C \times C) \times C}; (50,-15)*+{(B \times
B) \times B} }; {\ar_{\otimes \times 1} (50,-15)*+{(B \times B)
\times B}; (100,-15)*+{B \times B} }; {\ar^{a} (50,-15)*+{(B
\times B) \times B}; (50,-30)*+{B \times (B \times B)} }; {\ar^{F
\times (F \times F)} (0,-15)*+{C \times (C \times C)};
(50,-30)*+{B \times (B \times B)} }; {\ar_{F \times 1} (0,-30)*+{C
\times (B \times B)}; (50,-30)*+{B \times (B \times B)} };
{\ar@/^1.5pc/^{F \times 1} (0,-45)*+{C \times B}; (50,-50)*+{B
\times B} }; {\ar@/^2pc/^{1 \times \otimes} (50,-30)*+{B \times (B
\times B)}; (50,-50)*+{B \times B} }; {\ar@{=} (100,-30)*+{B};
(75,-50)*+{B} }; (90,-6)*{\scriptstyle \Downarrow \bs{\alpha}
\times 1}; (60,-6)*{=}; (30,-6)*{=}; (22,-15)*{=}; (22,-27)*{=};
(35,-38)*{=}; (26,-46)*{\scriptstyle \Downarrow 1 \times \bs{\alpha}^{\dot}};
(75,-30)*{\scriptstyle \Downarrow \bs{a}^{B}}
\endxy
\]

The left unit adjoint equivalence is given by the pasting below,
following the same conventions as above.
\[
\xy {\ar^{I^{B} \times 1} (0,0)*+{* \times C}; (40,0)*+{B \times
C} }; {\ar^{G \times 1} (40,0)*+{B \times C}; (80,0)*+{C \times C}
}; {\ar^{F \times F}  (80,0)*+{C \times C}; (80,-15)*{B \times B}
}; {\ar^{\otimes} (80,-15)*{B \times B}; (80,-30)*+{B} }; {\ar^{G}
(80,-30)*+{B}; (80,-45)*+{C} }; {\ar_{1 \times F} (40,0)*+{B
\times C}; (80,-15)*{B \times B} }; {\ar_{1 \times F} (0,0)*+{*
\times C}; (40,-15)*+{* \times B} }; {\ar_{I^{B} \times 1}
(40,-15)*+{* \times B}; (80,-15)*{B \times B} };
{\ar_{\textrm{proj}_{2}} (40,-15)*+{* \times B}; (80,-30)*+{B} };
{\ar_{1 \times G} (40,-15)*+{* \times B}; (40,-30)*+{* \times C}
}; {\ar_{\textrm{proj}_{2}}  (40,-30)*+{* \times C}; (80,-45)*+{C}
}; {\ar@/_2pc/_{1} (0,0)*+{* \times C}; (40,-30)*+{* \times C} };
(70,-6)*{\scriptstyle \Downarrow \bs{\alpha} \times 1};
(40,-8)*{=}; (65,-21)*{\scriptstyle \Downarrow \bs{l}^{B}}; (60,-31)*{=};
(20,-15)*{\scriptstyle \stackrel{\Leftarrow}{1 \times
\bs{\beta}^{\dot}}}
\endxy
\]
The right unit adjoint equivalence is constructed in an analogous
fashion.

In order to construct the remaining data, it is useful to compute
the components of the above transformations.  The component of
$a^{C}$ at $(x,y,z)$ is given as the following composite which we
write omitting all the identity components from adjoint
equivalences marked with an equal sign.
\[
\begin{array}{rcl}
G \big( FG(Fx \otimes Fy) \otimes Fz \big) & \stackrel{G \big( \alpha \otimes 1 \big)}{\longrightarrow} & G \big( (Fx \otimes Fy) \otimes Fz \big) \\
& \stackrel{Ga^{B}}{\longrightarrow} & G \big( Fx \otimes (Fy \otimes Fz) \big) \\
& \stackrel{G \big( 1 \otimes \alpha^{\dot}
\big)}{\longrightarrow} & G \big( Fx \otimes FG(Fy \otimes Fz)
\big)
\end{array}
\]
The component of $l^{C}$ at $x$ is given by the following
composite.
\[
G(FGI \otimes Fx) \stackrel{G(\alpha \otimes 1)}{\longrightarrow}
G(I \otimes Fx) \stackrel{Gl}{\rightarrow} GFx
\stackrel{\beta^{\dot}}{\rightarrow} x
\]
The component of $r^{C}$ at $x$ is given by the following
composite.
\[
G(Fx \otimes FGI) \stackrel{G(1 \otimes \alpha)}{\longrightarrow}
G(Fx \otimes I) \stackrel{Gr}{\rightarrow} GFx
\stackrel{\beta^{\dot}}{\rightarrow} x
\]

Now we define the invertible modifications $\pi, \mu, \lambda,
\rho$.  We explicitly define the unit modifications $\mu, \lambda,
\rho$, but only describe the construction of the modification
$\pi$ due to the size of the pasting diagram.

The unit modification $\mu$ is the pasting below (composed with
the unique coherence isomorphism $1 \circ 1 \circ 1 \cong 1$),
where all unmarked isomorphisms are naturality isomorphisms and
$\widetilde{\Phi}$ is the mate of $\Phi$.
\[
\xy {\ar^{\scriptstyle G(F \beta 1)} (0,0)*+{\scriptstyle
G(FxFy)}; (0,12)*+{\scriptstyle G(FGFxFy)} }; {\ar^{\scriptstyle
G(FG(r^{\dot})1)}  (0,12)*+{\scriptstyle G(FGFxFy)};
(5,24)*+{\scriptstyle G(FG(FxI)Fy)} }; {\ar^{\scriptstyle G(FG(1
\alpha^{\dot})1)} (5,24)*+{\scriptstyle G(FG(FxI)Fy)};
(15,34)*+{\scriptstyle G(FG(FxFGI)Fy)} }; {\ar^{\scriptstyle
G(\alpha 1)} (15,34)*+{\scriptstyle G(FG(FxFGI)Fy)};
(40,42)*+{\scriptstyle G((FxFGI)Fy)} }; {\ar^{\scriptstyle Ga}
(40,42)*+{\scriptstyle G((FxFGI)Fy)}; (70,42)*+{\scriptstyle
G(Fx(FGIFy))} }; {\ar^{\scriptstyle G(1 \alpha^{\dot})}
(70,42)*+{\scriptstyle G(Fx(FGIFy))}; (95,34)*+{\scriptstyle
G(FxFG(FGIFy))} }; {\ar^{\scriptstyle G(1FG(\alpha 1))}
(95,34)*+{\scriptstyle G(FxFG(FGIFy))}; (105,24)*+{\scriptstyle
G(FxFG(IFy))} }; {\ar^{\scriptstyle G(1FG(l))}
(105,24)*+{\scriptstyle G(FxFG(IFy))}; (110,12)*+{\scriptstyle
G(FxFGFy)} }; {\ar^{\scriptstyle G(1F \beta^{\dot})}
(110,12)*+{\scriptstyle G(FxFGFy)}; (110,0)*+{\scriptstyle
G(FxFy)} }; {\ar^{\scriptstyle G(\alpha 1)} (0,12)*+{\scriptstyle
G(FGFxFy)}; (35,0)*+{\scriptstyle G(FxFy)} }; {\ar^{\scriptstyle
G(r^{\dot} 1)} (35,0)*+{\scriptstyle G(FxFy)};
(34,23)*+{\scriptstyle G((FxI)Fy)} }; {\ar_{\scriptstyle G(\alpha
1)} (5,24)*+{\scriptstyle G(FG(FxI)Fy)}; (34,23)*+{\scriptstyle
G((FxI)Fy)} }; {\ar_{\scriptstyle G((1 \alpha^{\dot})1)}
(34,23)*+{\scriptstyle G((FxI)Fy)}; (40,42)*+{\scriptstyle
G((FxFGI)Fy)} }; {\ar_{\scriptstyle 1}  (0,0)*+{\scriptstyle
G(FxFy)}; (35,0)*+{\scriptstyle G(FxFy)} }; {\ar_{\scriptstyle Ga}
(34,23)*+{\scriptstyle G((FxI)Fy)}; (55,11)*+{\scriptstyle
G(Fx(IFy))} }; {\ar_{\scriptstyle G(1(\alpha^{\dot}1))}
(55,11)*+{\scriptstyle G(Fx(IFy))}; (70,42)*+{\scriptstyle
G(Fx(FGIFy))} }; {\ar^{\scriptstyle G(1 \alpha^{\dot})}
(55,11)*+{\scriptstyle G(Fx(IFy))}; (105,24)*+{\scriptstyle
G(FxFG(IFy))} }; {\ar_{\scriptstyle G(1l)} (55,11)*+{\scriptstyle
G(Fx(IFy))}; (70,0)*+{\scriptstyle G(FxFy)} }; {\ar_{\scriptstyle
1} (35,0)*+{\scriptstyle G(FxFy)}; (70,0)*+{\scriptstyle G(FxFy)}
}; {\ar^{\scriptstyle G(1 \alpha^{\dot})} (70,0)*+{\scriptstyle
G(FxFy)}; (110,12)*+{\scriptstyle G(FxFGFy)} }; {\ar_{\scriptstyle
1} (70,0)*+{\scriptstyle G(FxFy)}; (110,0)*+{\scriptstyle G(FxFy)}
}; (13,5)*{\scriptstyle \Downarrow \Phi}; (15,15)*{\scriptstyle
\cong}; (20,30)*{\scriptstyle \cong}; (55,30)*{\scriptstyle
\cong}; (80,30)*{\scriptstyle \cong}; (80,12)*{\scriptstyle
\cong}; (95,5)*{\scriptstyle \Downarrow \tilde{\Phi}}; (50,
7)*{\scriptstyle \Downarrow G \mu}
\endxy
\]

The unit modification $\lambda$ is the pasting below.
\[
\xy {\ar^{\scriptstyle G(FG(\alpha 1)1)} (0,0)*+{\scriptstyle
G(FG(FGIFx)Fy)}; (30,0)*+{\scriptstyle G(FG(IFx)Fy)} };
{\ar^{\scriptstyle G(FG(l)1)} (30,0)*+{\scriptstyle G(FG(IFx)Fy)};
(60,0)*+{\scriptstyle G(FGFxFy)} }; {\ar^{\scriptstyle G(F
\beta^{\dot}1)} (60,0)*+{\scriptstyle G(FGFxFy)};
(90,0)*+{\scriptstyle G(FxFy)} }; {\ar_{\scriptstyle G(\alpha 1)}
(0,0)*+{\scriptstyle G(FG(FGIFx)Fy)}; (0,-13)*+{\scriptstyle
G((FGIFx)Fy)} }; {\ar_{\scriptstyle Ga} (0,-13)*+{\scriptstyle
G((FGIFx)Fy)}; (0,-26)*+{\scriptstyle G(FGI(FxFy))} };
{\ar_{\scriptstyle G(1 \alpha^{\dot})} (0,-26)*+{\scriptstyle
G(FGI(FxFy))}; (0,-39)*+{\scriptstyle G(FGIFG(FxFy))} };
{\ar_{\scriptstyle G(\alpha 1)} (0,-39)*+{\scriptstyle
G(FGIFG(FxFy))}; (45,-39)*+{\scriptstyle G(IFG(FxFy))} };
{\ar@/_1pc/_{\scriptstyle Gl} (45,-39)*+{\scriptstyle
G(IFG(FxFy))}; (90,-26)*+{\scriptstyle GFG(FxFy)} };
{\ar_{\scriptstyle \beta^{\dot}} (90,-26)*+{\scriptstyle
GFG(FxFy)}; (90,0)*+{\scriptstyle G(FxFy)} }; {\ar_{\scriptstyle
G(\alpha 1)} (30,0)*+{\scriptstyle G(FG(IFx)Fy)};
(30,-13)*+{\scriptstyle G((IFx)Fy)} }; {\ar_{\scriptstyle
G((\alpha 1)1)} (0,-13)*+{\scriptstyle G((FGIFx)Fy)};
(30,-13)*+{\scriptstyle G((IFx)Fy)} }; {\ar_{\scriptstyle Ga}
(30,-13)*+{\scriptstyle G((IFx)Fy)}; (30,-26)*+{\scriptstyle
G(I(FxFy))} }; {\ar_{\scriptstyle G(\alpha(11))}
(0,-26)*+{\scriptstyle G(FGI(FxFy))}; (30,-26)*+{\scriptstyle
G(I(FxFy))} }; {\ar^{\scriptstyle G(1 \alpha^{\dot})}
(30,-26)*+{\scriptstyle G(I(FxFy))}; (45,-39)*+{\scriptstyle
G(IFG(FxFy))} }; {\ar_{\scriptstyle G(\alpha 1)}
(60,0)*+{\scriptstyle G(FGFxFy)}; (60,-19.5)*+{\scriptstyle
G(FxFy)} }; {\ar^{\scriptstyle G(l1)} (30,-13)*+{\scriptstyle
G((IFx)Fy)}; (60,-19.5)*+{\scriptstyle G(FxFy)} };
{\ar_{\scriptstyle Gl} (30,-26)*+{\scriptstyle G(I(FxFy))};
(60,-19.5)*+{\scriptstyle G(FxFy)} }; {\ar_{\scriptstyle 1}
(60,-19.5)*+{\scriptstyle G(FxFy)}; (90,0)*+{\scriptstyle G(FxFy)}
}; {\ar_{\scriptstyle G \alpha^{\dot}} (60,-19.5)*+{\scriptstyle
G(FxFy)}; (90,-26)*+{\scriptstyle GFG(FxFy)} };
(15,-6)*{\scriptstyle \cong}; (15,-19)*{\scriptstyle \cong};
(15,-32)*{\scriptstyle \cong}; (45,-6)*{\scriptstyle \cong}; (43,
-19)*{\scriptstyle \Downarrow G \lambda}; (60,-28)*{\scriptstyle
\cong}; (75,-5)*{\scriptstyle \Downarrow G(\tilde{\Phi}1)};
(78,-15)*{\scriptstyle \Downarrow \tilde{\Psi}}
\endxy
\]

For the unit modification $\rho:1 \otimes r^{\dot} \Rightarrow a
\circ r^{\dot}$, we give the following non-standard presentation.
Let $\tilde{\rho}$ denote the mate of $\rho^{-1}$ with source $(1
\otimes r) \circ a$ and target $r$.  We thus define
$\tilde{\rho}$ as the pasting below and leave it to the reader to
construct the usual modification $\rho$ as necessary.
\[
\xy {\ar^{\scriptstyle G(\alpha 1)} (0,0)*+{\scriptstyle
G(FG(FxFy)FGI)}; (30,0)*+{\scriptstyle G((FxFy)FGI)} };
{\ar^{\scriptstyle Ga}  (30,0)*+{\scriptstyle G((FxFy)FGI)};
(60,0)*+{\scriptstyle G(Fx(FyFGI))} }; {\ar^{\scriptstyle G(1
\alpha^{\dot})} (60,0)*+{\scriptstyle G(Fx(FyFGI))};
(90,0)*+{\scriptstyle G(FxFG(FyFGI))} }; {\ar^{\scriptstyle
G(1FG(1 \alpha))} (90,0)*+{\scriptstyle G(FxFG(FyFGI))};
(90,-13)*+{\scriptstyle G(FxFG(FyI))} }; {\ar^{\scriptstyle
G(1FGr)} (90,-13)*+{\scriptstyle G(FxFG(FyI))};
(90,-26)*+{\scriptstyle G(FxFGFy)} }; {\ar^{\scriptstyle G(1F
\beta^{\dot})} (90,-26)*+{\scriptstyle G(FxFGFy)};
(90,-39)*+{\scriptstyle G(FxFy)} }; {\ar_{\scriptstyle G(1
\alpha)} (0,0)*+{\scriptstyle G(FG(FxFy)FGI)};
(0,-26)*+{\scriptstyle G(FG(FxFy)I)} }; {\ar@/_1pc/_{\scriptstyle
Gr} (0,-26)*+{\scriptstyle G(FG(FxFy)I)}; (45,-39)*+{\scriptstyle
GFG(FxFy)} }; {\ar_{\scriptstyle \beta^{\dot}}
(45,-39)*+{\scriptstyle GFG(FxFy)}; (90,-39)*+{\scriptstyle
G(FxFy)} }; {\ar_{\scriptstyle G((11)\alpha)}
(30,0)*+{\scriptstyle G((FxFy)FGI)}; (30,-19.5)*+{\scriptstyle
G((FxFy)I)} }; {\ar^{\scriptstyle Ga} (30,-19.5)*+{\scriptstyle
G((FxFy)I)}; (60,-13)*+{\scriptstyle G(Fx(FyI))} };
{\ar^{\scriptstyle G(1(1 \alpha))} (60,0)*+{\scriptstyle
G(Fx(FyFGI))}; (60,-13)*+{\scriptstyle G(Fx(FyI))} };
{\ar_{\scriptstyle G(1 \alpha^{\dot})} (60,-13)*+{\scriptstyle
G(Fx(FyI))}; (90,-13)*+{\scriptstyle G(FxFG(FyI))} };
{\ar_{\scriptstyle Gr} (30,-19.5)*+{\scriptstyle G((FxFy)I)};
(60,-26)*+{\scriptstyle G(FxFy)} }; {\ar^{\scriptstyle G(\alpha
1)} (0,-26)*+{\scriptstyle G(FG(FxFy)I)};
(30,-19.5)*+{\scriptstyle G((FxFy)I)} }; {\ar^{\scriptstyle G(1
r)} (60,-13)*+{\scriptstyle G(Fx(FyI))}; (60,-26)*+{\scriptstyle
G(FxFy)} }; {\ar^{\scriptstyle G(1 \alpha^{\dot})}
(60,-26)*+{\scriptstyle G(FxFy)}; (90,-26)*+{\scriptstyle
G(FxFGFy)} }; {\ar^{\scriptstyle G \alpha} (45,-39)*+{\scriptstyle
GFG(FxFy)}; (60,-26)*+{\scriptstyle G(FxFy)} }; {\ar_{\scriptstyle
1} (60,-26)*+{\scriptstyle G(FxFy)}; (90,-39)*+{\scriptstyle
G(FxFy)} }; (15,-15)*{\scriptstyle \cong}; (45,-6)*{\scriptstyle
\cong}; (75,-6)*{\scriptstyle \cong}; (35,-26)*{\scriptstyle
\cong}; (48,-19.5)*{\scriptstyle \Downarrow G \tilde{\rho}};
(75,-19)*{\scriptstyle \cong}; (79,-30)*{\scriptstyle \Downarrow
G(1 \tilde{\Phi})}; (65,-33)*{\scriptstyle \Downarrow
\tilde{\Psi}}
\endxy
\]

We now describe the construction of the modification $\pi$ for the
monoidal structure on $C$.  This invertible modification is
constructed like the unit modifications by pasting together
\begin{itemize}
\item naturality isomorphisms for $\alpha, \alpha^{\dot}$, and the
associator $a$; \item a single counit coming from $\alpha
\dashv_{eq} \alpha^{\dot}$; and \item the cell $G(\pi^{B})$.
\end{itemize}
We leave it to the reader to construct the appropriate pasting from these cells.

There are now three monoidal bicategory axioms to check.  We leave
these to the reader as the diagrams are large but the computations
simple -- the associativity axiom follows by coherence,
naturality, and the associativity axiom in $B$, while both unit
axioms follow by coherence, naturality axioms, the corresponding
unit axioms in $B$, and the biadjoint biequivalence axioms.

Now we show that $F$ can be equipped with the structure of a
monoidal functor.  The adjoint equivalence
\[
\bs{\chi}: \otimes_{B} \circ (F \times F) \rightarrow F \circ
\otimes_{C}
\]
is $\bs{\alpha}^{\dot}$, and the adjoint equivalence
\[
\bs{\iota}:I_{B} \rightarrow F \circ I_{C}
\]
is $\bs{\alpha}_{I}^{\dot}$.  The invertible modification $\omega$
\[
\xy {\ar^{\chi \otimes 1} (0,0)*+{(FxFy)Fz}; (30,0)*+{F(xy)Fz} };
{\ar^{\chi}  (30,0)*+{F(xy)Fz}; (60,0)*+{F((xy)z)} }; {\ar^{Fa}
(60,0)*+{F((xy)z)}; (60,-20)*+{F(x(yz))} }; {\ar_{a}
(0,0)*+{(FxFy)Fz}; (0,-20)*+{Fx(FyFz)} }; {\ar_{1 \otimes \chi}
(0,-20)*+{Fx(FyFz)}; (30,-20)*+{FxF(yz)} }; {\ar_{\chi}
(30,-20)*+{FxF(yz)}; (60,-20)*+{F(x(yz))} }; (30,-10)*{\Downarrow
\omega}
\endxy
\]
is given by the pasting diagram below in which all the cells are
naturality isomorphisms or a counit for $FG(\alpha \otimes 1)
\dashv FG(\alpha^{\dot} \otimes 1)$.
\[
\xy
{\ar^{\scriptstyle \alpha^{\dot}1} (0,3)*+{\scriptstyle (FxFy)Fz}; (30,12)*+{\scriptstyle FG(FxFy)Fz}};
{\ar^{\scriptstyle \alpha^{\dot}} (30,12)*+{\scriptstyle FG(FxFy)Fz};
(60,3)*+{\scriptstyle FG(FG(FxFy)Fz)} };
{\ar^{\scriptstyle FG(\alpha 1)} (60,3)*+{\scriptstyle FG(FG(FxFy)Fz)}; (60,-12)*+{\scriptstyle FG((FxFy)Fz)} };
{\ar^{\scriptstyle FGa} (60,-12)*+{\scriptstyle FG((FxFy)Fz)}; (60,-24)*+{\scriptstyle FG(Fx(FyFz))} };
{\ar^{\scriptstyle FG(1 \alpha^{\dot})} (60,-24)*+{\scriptstyle FG(Fx(FyFz))};
(60,-36)*+{\scriptstyle FG(FxFG(FyFz))} };
{\ar_{\scriptstyle a} (0,3)*+{\scriptstyle (FxFy)Fz}; (0,-18)*+{\scriptstyle Fx(FyFz)} };
{\ar_{\scriptstyle 1 \alpha^{\dot}} (0,-18)*+{\scriptstyle Fx(FyFz)};
(0,-36)*+{\scriptstyle FxFG(FyFz)} };
{\ar_{\scriptstyle \alpha^{\dot}} (0,-36)*+{\scriptstyle FxFG(FyFz)}; (60,-36)*+{\scriptstyle FG(FxFG(FyFz))} };
{\ar_{\scriptstyle \alpha^{\dot}} (0,3)*+{\scriptstyle (FxFy)Fz}; (30,-6)*+{\scriptstyle FG((FxFy)Fz)} };
{\ar^{\scriptstyle FG(\alpha^{\dot} 1)} (30,-6)*+{\scriptstyle FG((FxFy)Fz)};
(60,3)*+{\scriptstyle FG(FG(FxFy)Fz)} };
{\ar_{\scriptstyle 1} (30,-6)*+{\scriptstyle FG((FxFy)Fz)};
(60,-12)*+{\scriptstyle FG((FxFy)Fz)} };
{\ar_{\scriptstyle \alpha^{\dot}}
(0,-18)*+{\scriptstyle Fx(FyFz)}; (60,-24)*+{\scriptstyle FG(Fx(FyFz))} }; (30,3)*{\scriptstyle \cong};
(50,-6)*{\scriptstyle \cong}; (30,-15)*{\scriptstyle \cong}; (30,-30)*{\scriptstyle \cong}
\endxy
\]
The invertible modification $\gamma$
\[
\xy {\ar^{\iota \otimes 1} (0,0)*+{IFx}; (20,0)*+{FIFx} };
{\ar^{\chi} (20,0)*+{FIFx}; (40,0)*+{F(Ix)} }; {\ar^{Fl}
(40,0)*+{F(Ix)}; (40,-20)*+{Fx} }; {\ar_{l} (0,0)*+{IFx};
(40,-20)*+{Fx} }; (25,-6)*{\Downarrow \gamma}
\endxy
\]
is the pasting diagram shown below where both cells are mates of
naturality isomorphisms.
\[
\xy {\ar^{\alpha^{\dot} 1} (0,0)*+{IFx}; (5,20)*+{FGIFx} };
{\ar^{\alpha^{\dot}} (5,20)*+{FGIFx}; (30,20)*+{FG(FGIFx)} };
{\ar^{FG(\alpha 1)} (30,20)*+{FG(FGIFx)}; (55,20)*+{FG(IFx)} };
{\ar^{FGl} (55,20)*+{FG(IFx)}; (80,20)*+{FGFx} }; {\ar^{\alpha}
(80,20)*+{FGFx}; (85,0)*+{Fx} }; {\ar_{l} (0,0)*+{IFx};
(85,0)*+{Fx} }; {\ar_{\alpha^{\dot}} (0,0)*+{IFx};
(55,20)*+{FG(IFx)} }; (22,12)*{\cong}; (55,6)*{\cong}
\endxy
\]
The invertible modification $\delta: \chi \circ (1 \otimes \iota)
\circ r^{\dot} \Rightarrow Fr^{\dot}$ is defined similarly.  There
are now two axioms to check to show that these data give a
monoidal functor between monoidal bicategories, and once again we
leave these simple albeit long computations to the reader.  The
associativity axiom follows from the transformation axioms for
$\alpha$ while the unit axiom requires using the biadjoint
biequivalence axioms.
\end{proof}

\remark
It would be possible at this point to prove that the functor $G$
chosen in the proof above can also be given the structure of a
monoidal functor.  For instance, the transformation $\chi$ is
given by the transformation $\beta$ used in constructing the
biadjoint biequivalence $F \dashv_{bieq} G$.  We could go further,
and even show that the entire biadjoint biequivalence $F
\dashv_{bieq} G$ can be lifted from $\textbf{Bicat}$ to
$\textbf{MonBicat}$, showing that the forgetful functor
\[
\textbf{MonBicat} \rightarrow \textbf{Bicat}
\]
is the tricategorical analogue of an isofibration.  Put another
way, the free monoidal bicategory construction is an example of a
``flexible 3-monad.''  The proofs of all of these statements
follow in exactly the same fashion as the construction of the
monoidal structure on $C$ and $F$ given above.
\endremark

\remark
We could also lift braided monoidal, sylleptic monoidal, or
symmetric monoidal structures along biequivalences in a similar
fashion.  For instance, if $B$ is braided with braiding
$R_{x,y}:xy \rightarrow yx$ then $C$ can be given a braided
structure with braiding
\[
F(GxGy) \stackrel{FR_{Gx,Gy}}{\longrightarrow} F(GyGx).
\]
In these cases, as above, the entire biadjoint biequivalence
could be lifted from $\textbf{Bicat}$ to the relevant tricategory
of monoidal bicategories of the kind considered.
\endremark

\section{Application:  Picard 2-categories}

This section will present an application of the main result to the
study of Picard 2-categories.  This is the analogue, for monoidal
bicategories, of the result of Baez and Lauda that the 2-category
of 2-groups (or Picard groupoids) is 2-equivalent to the
2-category of coherent 2-groups \cite{bl}.  It should be noted
that all of the results of this section remain true when we add
braided, sylleptic, or symmetric structures.

\definition A \textit{Picard 2-category} $X$ is a monoidal bicategory such that for every object $x$, there exists an object $y$ such that
\[
x \otimes y \simeq I \simeq y \otimes x.
\]
\enddefinition

\remark The reader should note that we call these Picard 2-categories even though the underlying object is a mere bicategory.  We have also not assumed that all the 1- and 2-cells are invertible, nor that the monoidal structure is braided, sylleptic, or symmetric.  All of these additional features (strictness, invertible higher cells, and symmetry) can be added as desired to produce the notion of Picard 2-category appropriate to a particular application.  Analogous results to those we present below can then be proven.
\endremark

\definition The tricategory $\mathbf{Pic2Cat}$ is the full sub-tricategory of $\mathbf{MonBicat}$ consisting of those monoidal bicategories which are Picard 2-categories.
\enddefinition

For the next definition, recall that every monoidal bicategory $X$ gives rise to a tricategory $\Sigma X$ with a single object $*$ and single hom-bicategory given by
\[
\Sigma X(*,*) = X;
\]
horizontal composition is then given by the tensor product, and all of the coherence constraints for the tricategory are given by those for the monoidal structure on $X$.  Thus a biadjoint biequivalence $x \dashv_{bieq} y$ between objects of a monoidal bicategory is defined to be a biadjoint biequivalence $x \dashv_{bieq} y$ in $\Sigma X$ where now $x,y$ are treated as 1-cells of the tricategory $\Sigma X$.

\definition A \textit{coherent Picard 2-category} $(X, \textrm{inv})$ is a monoidal bicategory $X$, a function $\textrm{inv}: \textrm{ob }X \rightarrow \textrm{ob }X$, and for each object $x$ a biadjoint biequivalence $x \dashv_{bieq} \textrm{inv}(x)$.
\enddefinition

\definition The tricategory $\mathbf{CohPic2Cat}$ has as its 0-cells coherent Picard 2-categories $X$, hom-bicategories defined as
\[
\mathbf{CohPic2Cat}\Big( (X, \textrm{inv}_{X}), (Y, \textrm{inv}_{Y}) \Big) = \mathbf{MonBicat}(X, Y),
\]
and all coherence constraints those inherited from the tricategory $\mathbf{MonBicat}$.
\enddefinition

\theorem
The underlying monoidal bicategory functor $U$ factors (as a strict functor between tricategories) through the inclusion of $\mathbf{Pic2Cat}$ into $\mathbf{MonBicat}$.
\[
\xy
{\ar^{U} (0,0)*+{\mathbf{CohPic2Cat}}; (40,0)*+{\mathbf{MonBicat}} };
{\ar@{.>}_{U'} (0,0)*+{\mathbf{CohPic2Cat}}; (20,-12)*+{\mathbf{Pic2Cat}} };
{\ar@{^{(}->} (20,-12)*+{\mathbf{Pic2Cat}}; (40,0)*+{\mathbf{MonBicat}} };
\endxy
\]
The strict functor $U': \mathbf{CohPic2Cat} \rightarrow \mathbf{Pic2Cat}$ is a triequivalence.
\endtheorem
\begin{proof}
The first statement is clear, as the underlying monoidal bicategory of a coherent Picard 2-category $(X, \textrm{inv})$ is obviously a Picard 2-category, and all of the higher dimensional structure involved in the definitions of these two tricategories agrees.  For the second statement, we must prove that $U'$ is locally a biequivalence and triessentially surjective.  Now $U'$ is the identity functor on hom-bicategories, so is locally a biequivalence.  To show that $U'$ is triessentially surjective, note that Theorem \ref{biadjbieqintricats} actually implies that $U'$ is surjective on objects as we can always choose a biadjoint biequivalence $x \dashv_{bieq} y$ for any object $x$ with the property that $x \otimes y \simeq I \simeq y \otimes x$.
\end{proof}

This theorem produces the most basic kind of equivalence between the theory of Picard 2-categories and its coherent version.  For the rest of this paper, we will sketch an improvement to this equivalence by explaining how one might go about proving that not only can Picard 2-categories be replaced by coherent ones, but also that monoidal functors can be replaced by ones that preserve the choice of inverses up to equivalence.

\definition  1.  Let $X$ be a bicategory.  Define the bicategory $X^{op}$ to be the one with
\begin{itemize}
\item the same objects as $X$,
\item $X^{op}(a,b) = X(b,a)$,
\item composition given by $g \circ^{op} f = f \circ g$, and
\item constraints given by $a_{h,g,f}^{op} = a_{f,g,h}^{-1}$, $l_{f}^{op} = r_{f}$, and $r_{f}^{op} = l_{f}$.
    \end{itemize}

2.  Let $Y$ be a monoidal bicategory.  Define the monoidal bicategory $Y^{rev}$ to be the one with
\begin{itemize}
\item underlying bicategory the same as $Y$,
\item $a \otimes^{rev} b = b \otimes a$,
\item $I^{rev} = I$,
\item all adjoint equivalences given by the opposites of the appropriate adjoint equivalences for the monoidal structure on $Y$, and
\item all invertible 2-cell data the same as that for $Y$.
\end{itemize}

3.  Let $Z$ be a monoidal bicategory.  Define $Z^{r}$ to be $(Z^{op})^{rev} = (Z^{rev})^{op}$.
\enddefinition

\proposition
Let $(X, \textrm{inv})$ be a coherent Picard 2-category.  The function on objects
\[
\textrm{inv}: \textrm{ob }X \rightarrow \textrm{ob }X
\]
extends to a monoidal functor of the same name,
\[
\textrm{inv}: X \rightarrow X^{r}.
\]
\endproposition
\begin{proof}
We have already defined $\textrm{inv}$ on objects, so now it is time to define it on 1- and 2-cells.  To do this, we must fix notation.  For an object $x$, the biadjoint biequivalence $x \dashv_{bieq} \textrm{inv}(x)$ consists of
\begin{itemize}
\item an adjoint equivalence $\mathbf{\epsilon}_{x} \dashv \mathbf{\epsilon}_{x}^{\dot}$ between the objects $x \otimes \textrm{inv}(x)$ and $I$,
\item an adjoint equivalence $\mathbf{\eta}_{x} \dashv \mathbf{\eta}_{x}^{\dot}$ between the objects $I$ and $\textrm{inv}(x) \otimes x$,
\item and invertible 2-cells $\Psi_{x}, \Phi_{x}$, satisfying the necessary axioms.
\end{itemize}
Thus for a 1-cell $f:x \rightarrow y$, we define $\textrm{inv}(f):\textrm{inv}(y) \rightarrow \textrm{inv}(x)$ as the following composite (ignoring associativity and unit constraints by coherence).
\[
\textrm{inv}(y) \stackrel{\eta_{x} 1}{\longrightarrow} \textrm{inv}(x) \otimes x \otimes     \textrm{inv}(y) \stackrel{1 f 1}{\longrightarrow} \textrm{inv}(x) \otimes y \otimes \textrm{inv}(y) \stackrel{1 \epsilon_{y}}{\longrightarrow} \textrm{inv}(x)
\]
We then define $\textrm{inv}(\alpha)$ to be $1*(1 \otimes \alpha \otimes 1) * 1$.

Next we must define structure constraints
\[
\begin{array}{c}
\textrm{inv}(g) \circ^{r} \textrm{inv}(f) \cong \textrm{inv}(g \circ f) \\
1_{x} \cong \textrm{inv}(1_{x})
\end{array}
\]
and check that these give a functor of bicategories.  Now $\textrm{inv}(g) \circ^{r} \textrm{inv}(f)$ in $X^{r}$ is defined to be the composite $\textrm{inv}(f) \circ \textrm{inv}(g)$, so we in fact require an isomorphism of the form
\[
\textrm{inv}(f) \circ \textrm{inv}(g) \cong \textrm{inv}(g \circ f).
\]
It is given by the pasting diagram below, in which we have written $\textrm{inv}(x)$ as $x^{-}$ and all the unmarked isomorphisms are functoriality of the tensor product.
\[
\xy
{\ar^{\scriptstyle \eta_{y}1} (0,0)*+{\scriptstyle z^{-}}; (30,0)*+{\scriptstyle y^{-}yz^{-}} };
{\ar^{\scriptstyle 1g1} (30,0)*+{\scriptstyle y^{-}yz^{-}}; (60,0)*+{\scriptstyle y^{-}zz^{-}} };
{\ar^{\scriptstyle 1 \epsilon_{z}} (60,0)*+{\scriptstyle y^{-}zz^{-}}; (90,0)*+{\scriptstyle y^{-}} };
{\ar^{\scriptstyle \eta_{x}1} (90,0)*+{\scriptstyle y^{-}}; (90,-14)*+{\scriptstyle x^{-}xy^{-}} };
{\ar^{\scriptstyle 1f1} (90,-14)*+{\scriptstyle x^{-}xy^{-}}; (90,-28)*+{\scriptstyle x^{-}yy^{-}} };
{\ar^{\scriptstyle 1 \epsilon_{y}} (90,-28)*+{\scriptstyle x^{-}yy^{-}}; (90,-42)*+{\scriptstyle x^{-}} };
{\ar_{\scriptstyle \eta_{x}1} (0,0)*+{\scriptstyle z^{-}}; (0,-14)*+{\scriptstyle x^{-}xz^{-}} };
{\ar_{\scriptstyle 11 \eta_{y}1} (0,-14)*+{\scriptstyle x^{-}xz^{-}}; (30,-14)*+{\scriptstyle x^{-}xy^{-}yz^{-}} };
{\ar_{\scriptstyle 111g1} (30,-14)*+{\scriptstyle x^{-}xy^{-}yz^{-}}; (60,-14)*+{\scriptstyle x^{-}xy^{-}zz^{-}} };
{\ar_{\scriptstyle 111\epsilon_{z}} (60,-14)*+{\scriptstyle x^{-}xy^{-}zz^{-}}; (90,-14)*+{\scriptstyle x^{-}xy^{-}} };
{\ar_{\scriptstyle 1f1} (0,-14)*+{\scriptstyle x^{-}xz^{-}}; (0,-28)*+{\scriptstyle x^{-}yz^{-}} };
{\ar^{\scriptstyle 11\eta_{y}1} (0,-28)*+{\scriptstyle x^{-}yz^{-}}; (30,-28)*+{\scriptstyle x^{-}yy^{-}yz^{-}} };
{\ar^{\scriptstyle 111g1} (30,-28)*+{\scriptstyle x^{-}yy^{-}yz^{-}}; (60,-28)*+{\scriptstyle x^{-}yy^{-}zz^{-}} };
{\ar^{\scriptstyle 111\epsilon_{z}} (60,-28)*+{\scriptstyle x^{-}yy^{-}zz^{-}}; (90,-28)*+{\scriptstyle x^{-}yy^{-}} };
{\ar@/_1.6pc/_{\scriptstyle 1} (0,-28)*+{\scriptstyle x^{-}yz^{-}}; (30,-42)*+{\scriptstyle x^{-}yz^{-}} };
{\ar_{\scriptstyle 1g1} (30,-42)*+{\scriptstyle x^{-}yz^{-}}; (60,-42)*+{\scriptstyle x^{-}zz^{-}} };
{\ar_{\scriptstyle 1\epsilon_{z}} (60,-42)*+{\scriptstyle x^{-}zz^{-}}; (90,-42)*+{\scriptstyle x^{-}} };
{\ar^{\scriptstyle 1\epsilon_{y}11} (30,-28)*+{\scriptstyle x^{-}yy^{-}yz^{-}}; (30,-42)*+{\scriptstyle x^{-}yz^{-}} };
{\ar^{\scriptstyle 1\epsilon_{y}11} (60,-28)*+{\scriptstyle x^{-}yy^{-}zz^{-}}; (60,-42)*+{\scriptstyle x^{-}zz^{-}} };
(45,-7)*{\scriptstyle \cong}; (45,-21)*{\scriptstyle \cong}; (45,-35)*{\scriptstyle \cong}; (75,-35)*{\scriptstyle \cong}; (15,-35)*{\scriptstyle \Downarrow 1 \Phi_{y} 1}
\endxy
\]
The isomorphism $1_{\textrm{inv}(x)} \cong \textrm{inv}(1_{x})$ is (modulo coherence) $\Phi_{x}^{-1}$.

There are now three axioms to check, one for associativity and two for units.  All three of these axioms follow by using the functoriality of the tensor product and then invoking coherence for monoidal bicategories (i.e., coherence for tricategories in the single-object case).  The unit axioms then require the equation $\Phi \circ \Phi^{-1} = 1$, while the associativity axiom only uses naturality.  Thus we have shown that inv is a functor of bicategories $X \rightarrow X^{op}$.  Now we turn to showing that it is monoidal.

The first step in showing that inv is a monoidal funtor $X \rightarrow X^{r}$ is to construct an adjoint equivalence
\[
\mathbf{\chi}: \otimes^{r} \circ (\textrm{inv} \times \textrm{inv}) \Rightarrow \textrm{inv} \circ \otimes.
\]
On components, this gives adjoint equivalences
\[
\textrm{inv}(x) \otimes^{r} \textrm{inv}(y) = \textrm{inv}(y) \otimes \textrm{inv}(x) \rightarrow \textrm{inv}(x \otimes y)
\]
in $X^{r}$, hence adjoint equivalences in $X$ with the source and target reversed.  We define the component $\chi_{x,y}$ below, and the rest of the adjoint equivalence will be defined in the obvious fashion.  We retain the same convention as above for writing $\textrm{inv}(x)$ as $x^{-}$.
\[
(xy)^{-} \stackrel{\eta_{y}1}{\longrightarrow} y^{-}y(xy)^{-} \stackrel{1\eta_{x}11}{\longrightarrow} y^{-}x^{-}xy(xy)^{-} \stackrel{11\epsilon_{xy}}{\longrightarrow} y^{-}x^{-}
\]

The second step in giving inv a monoidal structure is to construct an adjoint equivalence
\[
\mathbf{\iota}: I^{r} \Rightarrow \textrm{inv}(I)
\]
in $X^{r}$.  Since $I^{r} = I$, this is the obvious adjoint equivalence with left adjoint shown below.
\[
I^{-} \stackrel{l^{\dot}}{\longrightarrow} II^{-} \stackrel{\epsilon_{I}}{\longrightarrow} I
\]

The third step in giving inv a monoidal structure is to define three invertible modifications $\omega, \gamma, \delta$.  We leave it to the reader to write down the pasting diagrams as they are quite large, but we explain here which cells will be included in each.  In each case, there will be a large number of coherence cells from the monoidal bicategory structure, most of which will arise from the functoriality of the tensor product.  The other cells in each case are as follows:
\begin{itemize}
\item for $\omega$, the remaining cells are $\Phi_{xy}, \Psi_{\textrm{inv}(xyz)},$ and $\Phi^{-1}_{yz}$;
\item for $\gamma$, the remaining cells are $\Phi_{Ix}, \Phi_{I},$ and $\Psi_{\textrm{inv}(x)}$;
\item for $\delta$, the remaining cell is $\Phi_{I}$.
\end{itemize}
(The apparent assymetry in the definitions of $\gamma$ and $\delta$ is due to the fact that $\gamma$ has the cell $Fl$ in the source while $\delta$ has the cell $Fr^{\dot}$ in the target.)

Finally, there are two monoidal functor axioms to check.  Both follow from coherence for monoidal functors (a special case of coherence for functors of tricategories) and the biadjoint biequivalence axioms.
\end{proof}

\remark
One can go further to prove that the monoidal functor $\textrm{inv}$ is actually a monoidal biequivalence.  First we can equip $X^{r}$ with a canonical coherent Picard 2-category structure using the one given on $X$.  Then it is easy to show that $\textrm{inv}^{2}(x)$ is equivalent to $x$ in $X$ as both are weak inverses for $\textrm{inv}(x)$, from which it follows that $\textrm{inv}$ squares to the identity.
\endremark


\definition A \textit{coherent functor} $(F, c, u, v):(X, \textrm{inv}_{X}) \rightarrow (Y, \textrm{inv}_{Y})$ between coherent Picard 2-categories consists of the following data:
\begin{itemize}
\item a monoidal functor $F:X \rightarrow Y$,
\item an equivalence 1-cell $c_{x}:(Fx)^{-} \rightarrow F(x^{-})$ for each object $x \in X$, and
\item a pair of invertible 2-cells $u_{x}, v_{x}$ for each object $x \in X$ as displayed below.
    \[
    \xy
    {\ar^{\eta_{Fx}} (0,0)*+{I}; (0,15)*+{(Fx)^{-}Fx} };
    {\ar^{c_{x}1} (0,15)*+{(Fx)^{-}Fx}; (25,15)*+{F(x^{-})Fx} };
    {\ar^{\chi} (25,15)*+{F(x^{-})Fx}; (50,15)*+{F(x^{-}x)} };
    {\ar_{\iota} (0,0)*+{I}; (50,0)*+{FI} };
    {\ar_{F\eta_{x}} (50,0)*+{FI}; (50,15)*+{F(x^{-}x)} };
    (25,7)*{\Downarrow u_{x}}
    \endxy
    \]
    \[
    \xy
    {\ar^{\varepsilon_{Fx}} (0,0)*+{Fx(Fx)^{-}}; (0,15)*+{I} };
    {\ar^{\iota} (0,15)*+{I}; (50,15)*+{FI} };
    {\ar_{1c_{x}} (0,0)*+{Fx(Fx)^{-}}; (25,0)*+{FxF(x^{-})} };
    {\ar_{\chi} (25,0)*+{FxF(x^{-})}; (50,0)*+{F(xx^{-})} };
    {\ar_{F\varepsilon_{x}} (50,0)*+{F(xx^{-})}; (50,15)*+{FI} };
    (25,7)*{\Downarrow v_{x}}
    \endxy
    \]
    \end{itemize}
These are subject to the following axiom.
\[
\xy
{\ar^{\scriptstyle r} (20,0)*+{\scriptstyle FxI}; (0,15)*+{\scriptstyle Fx} };
{\ar^{\scriptstyle l^{\dot}} (0,15)*+{\scriptstyle Fx}; (20,30)*+{\scriptstyle IFx} };
{\ar^{\scriptstyle \iota 1} (20,30)*+{\scriptstyle IFx}; (60,30)*+{\scriptstyle FIFx} };
{\ar^{\scriptstyle \chi} (60,30)*+{\scriptstyle FIFx}; (100,30)*+{\scriptstyle F(Ix)} };
{\ar_{\scriptstyle 1 \iota} (20,0)*+{\scriptstyle FxI}; (60,0)*+{\scriptstyle FxFI} };
{\ar_{\scriptstyle \chi} (60,0)*+{\scriptstyle FxFI}; (100,0)*+{\scriptstyle F(xI)} };
{\ar_{\scriptstyle F(1 \eta_{x})} (100,0)*+{\scriptstyle F(xI)}; (100,15)*+{\scriptstyle F(xx^{-}x)} };
{\ar_{\scriptstyle F(\varepsilon_{x}1)} (100,15)*+{\scriptstyle F(xx^{-}x)}; (100,30)*+{\scriptstyle F(Ix)} };
{\ar_{\scriptstyle 1 \eta_{Fx}} (20,0)*+{\scriptstyle FxI}; (20,15)*+{\scriptstyle Fx(Fx)^{-}Fx} };
{\ar_{\scriptstyle \varepsilon_{Fx}1} (20,15)*+{\scriptstyle Fx(Fx)^{-}Fx}; (20,30)*+{\scriptstyle IFx} };
{\ar^{\scriptstyle 1c1} (20,15)*+{\scriptstyle Fx(Fx)^{-}Fx}; (45,15)*+{\scriptstyle FxF(x^{-})Fx} };
{\ar_{\scriptstyle 1 \chi} (45,15)*+{\scriptstyle FxF(x^{-})Fx}; (72.5,10)*+{\scriptstyle FxF(x^{-}x)} };
{\ar^{\scriptstyle \chi 1} (45,15)*+{\scriptstyle FxF(x^{-})Fx}; (72.5,20)*+{\scriptstyle F(xx^{-})Fx} };
{\ar^{\scriptstyle \chi} (72.5,20)*+{\scriptstyle F(xx^{-})Fx}; (100,15)*+{\scriptstyle F(xx^{-}x)} };
{\ar_{\scriptstyle \chi} (72.5,10)*+{\scriptstyle FxF(x^{-}x)}; (100,15)*+{\scriptstyle F(xx^{-}x)} };
{\ar_{\scriptstyle 1F\eta_{x}} (60,0)*+{\scriptstyle FxFI}; (72.5,10)*+{\scriptstyle FxF(x^{-}x)} };
{\ar_{\scriptstyle F\varepsilon_{x} 1} (72.5,20)*+{\scriptstyle F(xx^{-})Fx}; (60,30)*+{\scriptstyle FIFx} };
(8,15)*{\scriptstyle \stackrel{\Phi^{-1}_{Fx}}{\Rightarrow}}; (40,23)*{\scriptstyle \Downarrow v_{x}1}; (40,7)*{\scriptstyle \Downarrow 1u_{x}}; (84,23)*{\scriptstyle \cong}; (84,7)*{\scriptstyle \cong}; (72.5, 15)*{\scriptstyle \Downarrow \omega};
{\ar@{=} (60,-3)*{}; (60,-7)*{} };
{\ar^{\scriptstyle r} (20,-40)*+{\scriptstyle FxI}; (0,-25)*+{\scriptstyle Fx} };
{\ar^{\scriptstyle l^{\dot}} (0,-25)*+{\scriptstyle Fx}; (20,-10)*+{\scriptstyle IFx} };
{\ar^{\scriptstyle \iota 1} (20,-10)*+{\scriptstyle IFx}; (60,-10)*+{\scriptstyle FIFx} };
{\ar^{\scriptstyle \chi} (60,-10)*+{\scriptstyle FIFx}; (100,-10)*+{\scriptstyle F(Ix)} };
{\ar_{\scriptstyle 1 \iota} (20,-40)*+{\scriptstyle FxI}; (60,-40)*+{\scriptstyle FxFI} };
{\ar_{\scriptstyle \chi} (60,-40)*+{\scriptstyle FxFI}; (100,-40)*+{\scriptstyle F(xI)} };
{\ar_{\scriptstyle F(1 \eta_{x})} (100,-40)*+{\scriptstyle F(xI)}; (100,-25)*+{\scriptstyle F(xx^{-}x)} };
{\ar_{\scriptstyle F(\varepsilon_{x}1)} (100,-25)*+{\scriptstyle F(xx^{-}x)}; (100,-10)*+{\scriptstyle F(Ix)} };
{\ar^{\scriptstyle r} (20,-40)*+{\scriptstyle FxI}; (60,-25)*+{\scriptstyle Fx} };
{\ar^{\scriptstyle l} (20,-10)*+{\scriptstyle IFx}; (60,-25)*+{\scriptstyle Fx} };
{\ar^{\scriptstyle Fr} (100,-40)*+{\scriptstyle F(xI)}; (60,-25)*+{\scriptstyle Fx} };
{\ar^{\scriptstyle Fl^{\dot}} (60,-25)*+{\scriptstyle Fx}; (100,-10)*+{\scriptstyle F(xI)} };
(20,-25)*{\scriptstyle \cong}; (60,-35)*{\scriptstyle \Downarrow \delta}; (60,-15)*+{\scriptstyle \Downarrow \gamma}; (80,-25)*{\scriptstyle \stackrel{\Rightarrow}{F\Phi^{-1}_{x}}}
\endxy
\]
\enddefinition

\remark
We could have structured this definition in a slightly different fashion in a variety of ways.  First, the equivalence 1-cells $c_{x}$ could have been the components of a transformation $c$ as shown here.
    \[
\xy
{\ar^{F} (0,0)*+{X}; (30,0)*+{Y} };
{\ar^{\textrm{inv}} (30,0)*+{Y}; (30,-15)*+{Y^{r}} };
{\ar_{\textrm{inv}} (0,0)*+{X}; (0,-15)*+{X^{r}} };
{\ar_{F^{r}} (0,-15)*+{X^{r}}; (30,-15)*+{Y^{r}} };
(15,-7.5)*{\Downarrow c}
\endxy
\]

 Second, we could have asked that the invertible 2-cells $u_{x}, v_{x}$ could have been the components of a pair of invertible modifications.  To express the axioms above as diagrams of modifications would have required that the $\eta_{x}, \varepsilon_{x}$ be the components of transformations, which in turn would require that the functor $\textrm{inv}$ be covariant instead of contravariant.  Thus we would have to restrict attention to those Picard 2-categories in which every 1-cell is an equivalence; in fact, we would need every 1-cell to come as part of a specified adjoint equivalence in order to prescribe $\textrm{inv}$ as a covariant \textit{functor}.

 Third, we could have required a third axiom about how $c, u, v$ interact with $\Psi$.  This axiom follows from the first axiom by using the biadjoint biequivalence axioms.  In addition, the pastings involved are larger than the one for the axiom above as they involve two different uses of $c$ instead of just one, so it requires additional naturality squares.
\endremark

\theorem
Let $(X, \textrm{inv}_{X}), (Y, \textrm{inv}_{Y})$ be coherent Picard 2-categories, and let $F:X \rightarrow Y$ be a monoidal functor.  Then $F$ underlies a coherent functor \[
(F,c, u, v): (X, \textrm{inv}_{X}) \rightarrow (Y, \textrm{inv}_{Y}).
\]
\endtheorem
\begin{proof}
The 1-cell $c_{x}:F(x^{-}) \rightarrow Fx^{-}$ is given by the following composite where the last arrow is given by a composite of coherence cells and is thus unique up to unique isomorphism by coherence for functors.
\[
F(x^{-}) \stackrel{\eta_{Fx}1}{\longrightarrow} Fx^{-}FxF(x^{-}) \stackrel{1 \chi}{\longrightarrow} Fx^{-}F(xx^{-}) \stackrel{1 F\varepsilon}{\longrightarrow} Fx^{-}FI \longrightarrow Fx^{-}
\]
It is immediate that $c_{x}$ is an equivalence 1-cell.

We must now construct the invertible 2-cells $u_{x}, v_{x}$ and check the two axioms.  The cell $u_{x}$ is a pasting of coherence cells from both $Y$ and the functor $F$, together with $\Phi_{Fx}$.  The cell $v_{x}$ is constructed similarly out of coherence cells and $F\Phi_{x}^{-1}$.  The two axioms are straightforward diagram chases.
\end{proof}

\corollary
Let $X$ be a Picard 2-category.  Then $X$ has a coherent structure which is unique in the following sense:  if $(X, \textrm{inv})$ and $(X, \textrm{inv}')$ are two coherent structures on $X$, then the identity functor on $X$ lifts to a coherent functor
\[
(1,c, u, v): (X, \textrm{inv}) \rightarrow (X, \textrm{inv}').
\]
\endcorollary

\remark
We leave it to the reader to define coherent transformations and modifications.  A coherent transformation will involve additional data, while a coherent modification will only involve a new axiom.  Defined correctly, it is then possible to prove that the forgetful functor from the tricategory in which all cells are coherent to the tricategory $\mathbf{Pic2Cat}$ is a triequivalence.  This shows that a coherent structure on a given Picard 2-category is unique in the strongest possible sense.
\endremark

\providecommand{\bysame}{\leavevmode\hbox
to3em{\hrulefill}\thinspace}
\providecommand{\MR}{\relax\ifhmode\unskip\space\fi MR }
\providecommand{\MRhref}[2]{%
  \href{http://www.ams.org/mathscinet-getitem?mr=#1}{#2}
} \providecommand{\href}[2]{#2}

\end{document}